\documentclass[onecolumn]{autart}
\usepackage{mathrsfs}    

\usepackage{graphicx}           
\usepackage{dcolumn}            
\usepackage{bm}                 

\usepackage{graphics}           
\usepackage{epsfig}             
\usepackage{times}              
\usepackage[fleqn]{amsmath}
\usepackage{amssymb}
\usepackage{extarrows}
\usepackage{amsfonts}
\usepackage{fancyhdr}
\usepackage{color}
\usepackage{subfigure}
\usepackage{enumerate}

\allowdisplaybreaks[4] 

\newtheorem{theorem}{Theorem}
\newtheorem{lemma}[theorem]{Lemma}

\newtheorem{proposition}[theorem]{Proposition}

\newtheorem{remark}{Remark}

\begin{document}

\begin{frontmatter}

\title{Regularity in the two-phase free boundary problems under non-standard
     growth conditions} 

\author{Jun Zheng}\ead{zhengjun2014@aliyun.com,zhengjun@swjtu.edu.cn}

\address{School of Mathematics, Southwest Jiaotong University,
        Chengdu 611756, China}
\begin{abstract}                           
In this paper, we prove several regularity results for the heterogeneous, two-phase free boundary problems $\mathcal {J}_{\gamma}(u)=\int_{\Omega}\big(f(x,\nabla
u)+\lambda_{+}(u^{+})^{\gamma}+\lambda_{-}(u^{-})^{\gamma}+gu\big)\text{d}x\rightarrow \text{min}$
under non-standard growth conditions. Included in such problems are
heterogeneous jets and cavities of Prandtl-Batchelor type with
$\gamma=0$, chemical reaction problems with $0<\gamma<1$, and obstacle
type problems with $\gamma=1$. Our results hold not only in
the degenerate case of $p> 2$ for $p-$Laplace equations, but
also in the singular
case of $1<p<2$, which are extensions of \cite{LdT}.
\end{abstract}
\end{frontmatter}
\date{}
 \section{Introduction}
Let $\Omega$ be a bounded open set
in $\mathbb{R}^{n}(n\geq 2)$, and $g\in L^{q}(\Omega),\psi \in W^{1,p}(\Omega)\cap L^{\infty}(\Omega)$ with $ \psi ^+=\max\{\pm \psi,0\}\neq 0$ and $p\geq 2,q\geq n$. In \cite{LdT}, Leit\~{a}o, de Queiroz and Teixeira provided
a complete description of the sharp regularity of minimizers to the heterogeneous, two-phase free boundary problems
\begin{align}\label{+1}
\mathscr{J}_{\gamma}(u)=\int_{\Omega}\big(|\nabla
u|^p+F_{\gamma}(u)+gu\big)\text{d}x \rightarrow \min,
\end{align}
over the set $\{u \in W^{1,p}(\Omega):u-\psi\in W^{1,p}_0(\Omega)\}$, where
\begin{align*}
F_{\gamma}(u)=\lambda_{+}(u^{+})^{\gamma}+\lambda_{-}(u^{-})^{\gamma},
\end{align*}
$\gamma \in [0,1]$ is a parameter, $0<\lambda_{-}<\lambda_{+}<+\infty$, and by convention,
\begin{align*}
F_{0}(u)=\lambda_{+}\chi_{\{u>0\}}+\lambda_{-}\chi_{\{u\leq 0\}}.
\end{align*}
The lower limiting case, i.e., $\gamma=0$, relates to jets
and cavities problems. The upper case, i.e., $\gamma=1$, relates to
obstacle type problems. The intermediary problem, i.e., $0<\gamma<1$, can
be used to model the density of certain chemical specie, in reaction
with a porous catalyst pellet. The authors established local
$C^{1,\alpha}-$ and Log-Lipschitz regularities for minimizers of the
functional $\mathscr{J}_{\gamma}$ when $\gamma\in(0,1],q>n$ and $\gamma=0,q=n$
in \eqref{+1} respectively, see \cite{LdT}.

Problem \eqref{+1} was extended to a large class of the following heterogeneous, two-phase free boundary problems in \cite{ZFZ,ZZZ}
\begin{align*}
\int_{\Omega}(A(|\nabla
u|)+F_{\gamma}(u)+gu)\text{d}x\rightarrow \min,
\end{align*}
 over the set $\{u\in W^{1,A}(\Omega):u-\psi\in W^{1,A}_{0}(\Omega)\}$,
  for given functions $g\in L^{\infty}(\Omega)$ and $\psi \in W^{1,A}(\Omega) L^{\infty}(\Omega)$ with $\psi ^+\neq 0$, where
  $W^{1,A}(\Omega)$ is the class of weakly differentiable functions
  with $\int_{\Omega}A(|\nabla
u|)\text{d}x<\infty$. Under Lieberman's condition on $A$, which allows for a
 different behavior at $0$ and at $\infty$, local Log-Lipschitz continuity and local $C^{1,\alpha}-$regularity of minimizers have been obtained for $\gamma=0$ and $\gamma\in(0,1]$ respectively in
the setting of Orlicz spaces, see \cite{ZFZ,ZZZ}.

The aim of this paper is to study the heterogeneous, two-phase free boundary problems \begin{align}\label{1.1}
\mathcal {J}_{\gamma}(u)=\int_{\Omega}\big(f(x,\nabla
u)+F_{\gamma}(u)+gu\big)\text{d}x\rightarrow \text{min},
\end{align}
over the set $\{u \in W^{1,p(\cdot)}(\Omega):u-\psi\in W^{1,p(\cdot)}_0(\Omega)\}$ in the framework of Sobolev spaces with variable exponents, where $f:\Omega \times
\mathbb{R}^{n}\rightarrow \mathbb{R}$ is a Carath\'{e}odory function having a form:
\begin{align}
L^{-1}|z|^{p(x)}\leq f(x,z)\leq L(1+|z|^{p(x)}),\notag
\end{align}
for all $x\in \Omega,z\in \mathbb{R}^{n}$, with $p:\Omega\rightarrow
(1,+\infty)$ a continuous function and $L\geq 1$ a constant. We
  establish local Log-Lipschitz continuity and local $C^{1,\alpha}-$regularity for minimizers of
  $\mathcal {J}_{\gamma}$ with $\gamma=0$, and $\gamma\in(0,1]$ respectively.

To the knowledge of the author, the present paper seems to be a first regularity result for the heterogeneous, two-phase free boundary problems \eqref{1.1} with $p(x)-$growth. It should be mentioned that a large class of functionals and identical obstacle problems under non-standard growth conditions
 have been
 studied in \cite{E2,EH1,EH2,EH3,AM}, which provide the reference estimates, and suitable localization
and freezing techniques, etc., to treat the nonstandard growth exponents in the functional governed by \eqref{1.1}. The results obtained in this paper are not only extensions of one-phase obstacle problems under non-standard
  growth conditions (see, e.g., \cite{EH2,EH3}), but also a supplement of the
  degenerate two-phase free boundary problems studied in \cite{LdT}, since
  our results contain the singular case of $1<p<2$.


   The rest of this paper is organized as follows. In section 2, we present some basic notations,
   definitions, assumptions, and the main results obtained in this paper,
   including existence and $L^{\infty}$-boundedness results (Theorem 2.1), and local H\"{o}lder, $C^{1,\alpha}-$ and Log-Lipschitz regularities
    of minimizers (Theorem 2.2 - 2.4). In Section 3, we carry out the existence and $L^{\infty}-$boundedness for
    minimizers of the functional $\mathcal {J}_{\gamma}
    (\gamma\in[0,1])$. In Section 4, we establish the higher integrability for minimizers of
     the functional $\mathcal {J}_{\gamma}
    (\gamma\in[0,1])$. In Section 5, we address local $C^{0,\alpha}-$regularity for minimizers
     of the functional having a form $\int_\Omega \big(h(\nabla u)+F_{\gamma}(u)+gu\big) \text{d}x$
     with $\gamma\in[0,1]$ (Theorem 2.2), where $h$ satisfies certain non-standard growth conditions. In Section 6, we prove local $C^{0,\alpha}-$regularity for minimizers of the
       functional $\mathcal {J}_{\gamma}(\gamma\in[0,1])$ (Theorem
       2.3). In Section 7 and 8, we establish local $C^{1,\alpha}-$regularity for minimizers of the
       functional $\mathcal {J}_{\gamma}(\gamma\in(0,1])$ and
       local Log-Lipschitz continuity for minimizer of
$\mathcal {J}_{0}$, (Theorem
       2.4) respectively.

\section{Preliminaries and Statements}
In this paper, $\Omega$ will denote an open bounded domain in
$\mathbb{R}^{n} (n\geq 2)$ and $B_{R}(x)$ the open ball
$\{y\in\mathbb{R}^{n}:|x-y|<R\}$ with centre $x\in \mathbb{R}^{n}$. If $u$ is an integrable function
defined on $B_{R}(x)$, we will set
$
(u)_{x,R}=\frac{\int_{B_{R}(x)}u(x)\text{d}x}{|B_{R}(x)|},\notag
$
where $|B_{R}(x)|$ is the Lebesgue measure of $B_{R}(x)$. Without
confusion, we will write $B_{R}$ and $(u)_{R}$ instead of $B_{R}(x)$
and $(u)_{x,R}$ respectively. We may write $C$ or $c$ as
a constant that may be different from each other, but independent of
$\gamma$.

  Let $p:\Omega\rightarrow
(1,+\infty)$ be a continuous function. The variable exponent Lebesgue space $L^{p(\cdot)}(\Omega)$ is
  defined by
$
L^{p(\cdot)}(\Omega)=\{u|\ u:\Omega\rightarrow\mathbb{R}\ \text{is\
measurable},\ \int_{\Omega}|u|^{p(x)}\text{d}x<+\infty\},
$
with the norm
$
\|u\|_{L^{p(\cdot)}(\Omega)}=\inf\{\lambda>0;\
\int_{\Omega}\left|\frac{u}{\lambda}\right|^{p(x)}\text{d}x\leq
1\}.
$
The variable exponent Sobolev space $W^{1,p(\cdot)}(\Omega)$ is
  defined by
$
W^{1,p(\cdot)}(\Omega)=\{u\in L^{p(\cdot)}(\Omega):|\nabla u|\in
L^{p(\cdot)}(\Omega),\notag
$
with the norm
$
\|u\|_{W^{1,p(\cdot)}(\Omega)}=\|u\|_{L^{p(\cdot)}(\Omega)}+\|\nabla
u\|_{L^{p(\cdot)}(\Omega)}.\notag
$
Define $W^{1,p(\cdot)}_{0}(\Omega)$ as the closure of
$C^{\infty}_{0}(\Omega)$ in $W^{1,p(\cdot)}(\Omega)$. We point out
that, if $\Omega$ is bounded and $p(\cdot)$ satisfies \eqref{2.6}, then
the spaces $L^{p(\cdot)}(\Omega),W^{1,p(\cdot)}(\Omega)$ and
  $W^{1,p(\cdot)}_{0}(\Omega)$ are  all separable and reflexive
  Banach spaces. $\|\nabla
u\|_{L^{p(\cdot)}(\Omega)}$ is an equivalent norm on
$W^{1,p(\cdot)}_{0}(\Omega)$. We refer to \cite{KR,FSZ,FZ} for more details of
  the space $W^{1,p(\cdot)}(\Omega)$.

  In this paper, we consider the following growth, ellipticity and continuity
 conditions:
\begin{align}
&f:\Omega\times\mathbb{R}^{n}\rightarrow \mathbb{R},\
{f(x,z)\ \text{is} \ C^2-\text{continuous\ in}\ x \ \text{and}\ z, \ \text{and\ convex\ in}\ z \text{\ for\ every}\ x},\label{2.1}\\
&L^{-1}(\mu^{2}+|z|^{2})^{\frac{p(x)}{2}}\leq f(x,z)\leq
L(\mu^{2}+|z|^{2})^{\frac{p(x)}{2}},\label{2.2}\\
&|f(x,z)-f(x_{0},z)|\leq
Lw(|x-x_{0}|)[(\mu^{2}+|z|^{2})^{\frac{p(x)}{2}}
+(\mu^{2}+|z|^{2})^{\frac{p(x_{0})}{2}}][1+\log(\mu^{2}+|z|^{2})],\label{2.3}
\end{align}
for all $z\in\mathbb{R}^{n},x$ and $x_{0}\in \Omega$, where $L\geq
1,\mu\in[0,1]$, $w:\mathbb{R}^{+}\rightarrow\mathbb{R}^{+}$ is
a nondecreasing continuous function, vanishing at zero, which
represents the modulus of $p$,
\begin{align}\label{2.4}
|p(x)-p(y)|\leq w(|x-y|)\ \ \text{for\ all}\ x,y\in
\overline{\Omega},
\end{align}
 and satisfying
$
\limsup\limits_{R\rightarrow0}\omega(R)\log\big(\frac{1}{R}\big)<+\infty,\notag
$
thus without loss of generality, assume that
\begin{align}\label{2.5}
\omega(R)\leq L|\log R|^{-1},
\end{align}
for all $R<1$. Moreover, we assume that
\begin{align}\label{2.6}
1<p_{-}=\inf_{x\in\Omega}p(x) \leq
p(x)\leq \sup\limits_{x\in\Omega}p(x)=p_{+} <+\infty\ \ \text{for\ all}\
x\in\Omega.
\end{align}
  Let $q:\Omega\rightarrow(1,+\infty)$ be a continuous function
 fulfilling the conditions of the type \eqref{2.4} and
 \eqref{2.5}. We always make the following assumptions on $p(\cdot)$
 and $q(\cdot)$:
\begin{align}\label{2.7}
 \frac{1}{p_{-}}-\frac{1}{p_{+}}<\frac{1}{n},\ \
 q(x)\geq q_{-}\ \text{for\ all}\ x\in \Omega,\ \
  q_{-}>\left\{
\begin{array}{l}
\frac{1}{p_{-}-1}\frac{1}{\frac{1}{n}-\frac{1}{p_{-}}+\frac{1}{p_{+}}}>n,\ \text{if}\ p_{-}<2, \\
\frac{1}{\frac{1}{n}-\frac{1}{p_{-}}+\frac{1}{p_{+}}}\geq n,\ \text{if}\ p_{-}\geq 2.%
\end{array}%
\right.
 \end{align}
  Given $\psi\in W^{1,p(\cdot)}(\Omega)\cap L^{\infty}(\Omega)$ and $g\in
 L^{q(\cdot)}(\Omega)$, let $\mathcal {K}=\{u \in W^{1,p(\cdot)}(\Omega) ;u-\psi\in W^{1,p(\cdot)}_0(\Omega)\}$.  
We say that a function $u\in \mathcal {K}$ is a minimizer of the functional
   $\mathcal {J}_{\gamma}(u)$ governed by \eqref{1.1} if
    $\mathcal {J}_{\gamma}(u)\leq \mathcal {J}_{\gamma}(v)$ for all
    $v\in \mathcal {K}$.

The first result obtained in this paper concerns with the existence and $L^\infty-$boundedness of minimizers of $\mathcal {J}_{\gamma}(u)$ governed by \eqref{1.1}.
 \begin{theorem}\label{Theorem 2.1}
Under assumptions \eqref{2.1}-\eqref{2.7}, for each
$0\leq \gamma\leq 1$, there exists a minimizer $u_{\gamma}\in
\mathcal {K}$ of the functional $\mathcal {J}_{\gamma}(u)$ governed
by \eqref{1.1}. Furthermore, $u_{\gamma}$ is bounded. More precisely,
\begin{align}
\|u_{\gamma}\|_{L^{\infty}(\Omega)}\leq
C(n,L,q_{-},p_{\pm},\lambda_{\pm},\Omega,\|\psi\|_{L^{\infty}(\partial\Omega)},\|g\|_{L^{q(\cdot)}(\Omega)}).\notag
\end{align}
\end{theorem}
  Now let \begin{align}\label{2.9}
\mathcal {H}_{\gamma}(u)=\int_{\Omega}\big(h(\nabla
u)+F_{\gamma}(u)+gu\big)\text{d}x,
\end{align}
where $h:\mathbb{R}^{n}\rightarrow\mathbb{R}$ is a $C^2$-\text{continuous} \ \text{and\ convex\ } function satisfying for all $z\in\mathbb{R}^{n}$,
 \begin{align}
L^{-1}(\mu^{2}+|z|^{2})^{\frac{p(x)}{2}}\leq h(z)\leq
L(\mu^{2}+|z|^{2})^{\frac{p(x)}{2}}.\label{2.2'}
\end{align}
  We present then the regularity properties of minimizers of
the functionals $\mathcal {H}_{\gamma}$ and $\mathcal
{J}_{\gamma}$.
\begin{theorem}\label{Theorem 2.2}
Assume that \eqref{2.2'} and \eqref{2.4}-\eqref{2.7} hold. If $u_{\gamma}\in \mathcal
{K}$ is a minimizer of the functional $\mathcal {H}_{\gamma} (\gamma\in[0,1])$ governed
by \eqref{2.9}, then $u_{\gamma}\in C^{0,\alpha}_{loc}(\Omega)$ for some
$\alpha\in(0,1)$.
\end{theorem}
\begin{theorem}\label{Theorem 2.3}
Assume that  \eqref{2.1}-\eqref{2.7} hold. If
$u_{\gamma}\in \mathcal {K}$ is a minimizer of
the functional $\mathcal {J}_{\gamma}(\gamma\in[0,1])$ governed by \eqref{1.1}, then
$u_{\gamma}\in C^{0,\alpha}_{loc}(\Omega)$ for some
$\alpha\in(0,1)$.
\end{theorem}
\begin{theorem}\label{Theorem 2.4}
Assume that \eqref{2.1}-\eqref{2.7} hold, and assume further that
$
\omega(R)\leq LR^{\varsigma}
$
for some $\varsigma>\frac{n}{q_{-}}\frac{p_{-}}{p_{-}-1}$ and all
$R\leq 1$. The following
statements hold true:
\begin{enumerate}[(i)]
  \item
For each $\gamma\in(0,1]$, every minimizer $u_{\gamma}$ of the
functional $\mathcal {J}_{\gamma}$ governed by \eqref{1.1} is
$C^{1,\alpha}_{loc}-$continuous for some $\alpha\in(0,1)$.
  \item For each $\gamma=0$, every minimizer $u_{0}$ of the functional $\mathcal {J}_{0}$ governed by \eqref{1.1} is
locally Log-Lipschitz continuous in $\Omega$, and therefore is
$C^{0,\alpha}_{loc}-$continuous for any $\alpha\in(0,1)$.
\end{enumerate}
\end{theorem}
\section{Existence and $L^{\infty}$-boundedness of minimizers}
In this section, we establish the existence and $L^{\infty}$-boundedness
for minimizers of the functional $\mathcal{J}_{\gamma}(\gamma\in[0,1])$.
\begin{pf*}{Proof of Theorem~\ref{Theorem 2.1}}Firstly we consider the existence
of a minimizer of the functional $\mathcal {J}_{\gamma}$. Let $I_{0}=\min\{\mathcal {J}_{\gamma}(u):u\in
\mathcal {K}\}$. Initially we claim that
$I_{0}>-\infty$. Indeed, for any $u\in\mathcal {K}$, by Poincar\'{e}'s inequality
there exists a positive constant $C=C(n,p_{\pm},\Omega)$ such that
\begin{align}\label{3.1}
\|u\|_{L^{p(\cdot)}(\Omega)}\leq &~\|u-\psi\|_{L^{p(\cdot)}(\Omega)}
+\|\psi\|_{L^{p(\cdot)}(\Omega)}\notag\\
\leq &~ C\|\nabla
u-\nabla\psi\|_{L^{p(\cdot)}(\Omega)}+\|\psi\|_{L^{p(\cdot)}(\Omega)}\notag\\
\leq &~C(\|\nabla
u\|_{L^{p(\cdot)}(\Omega)}+\|\nabla\psi\|_{L^{p(\cdot)}(\Omega)}+\|\psi\|_{L^{p(\cdot)}(\Omega)}),
\end{align}
which implies
\begin{align}\label{3.2}
\|\nabla u\|^{p_{-}}_{L^{p(\cdot)}(\Omega)}\geq C_{1}\|
u\|^{p_{-}}_{L^{p(\cdot)}(\Omega)}-\|
\psi\|^{p_{-}}_{L^{p(\cdot)}(\Omega)}-\|\nabla\psi\|^{p_{-}}_{L^{p(\cdot)}(\Omega)},
\end{align}
and
\begin{align}\label{3.3}
\|\nabla u\|^{p_{+}}_{L^{p(\cdot)}(\Omega)}\geq C_{2}\|
u\|^{p_{+}}_{L^{p(\cdot)}(\Omega)}-\|
\psi\|^{p_{+}}_{L^{p(\cdot)}(\Omega)}-\|\nabla\psi\|^{p_{+}}_{L^{p(\cdot)}(\Omega)},
\end{align}
where $C_{1},\ C_{2}$ are positive constants depending only on $n,\
p_{\pm},\ \Omega$.

 Due to $q(x)\geq q_{-}$, we deduce by \eqref{2.7} and H\"{o}lder's
inequality that
\begin{align}
\bigg|\int_{\Omega}gu\text{d}x\bigg|\leq &
C_{3}(p_{+},p_{-})\|g\|_{L^{\frac{p(\cdot)}{p(\cdot)-1}}(\Omega)}\|u\|_{L^{p(\cdot)}(\Omega)}\notag\\
\leq &~
C_{4}(p_{+},p_{-})\|g\|_{L^{q(\cdot)}(\Omega)}\|1\|_{L^{\frac{1}{1-\frac{1}{p(\cdot)}
-\frac{1}{q(\cdot)}}}(\Omega)}\|u\|_{L^{p(\cdot)}(\Omega)}\notag\\
\leq &~ C_{4}(p_{+},p_{-})\bigg(1+|\Omega|^{1-\frac{1}{p_{-}}
-\frac{1}{q_{-}}}\bigg)\|g\|_{L^{q(\cdot)}(\Omega)}\|u\|_{L^{p(\cdot)}(\Omega)}\label{3.4}\\
\leq &~\begin{cases} \varepsilon\|u\|^{p_{-}}_{L^{p(\cdot)}(\Omega)}+
C_{5}(\varepsilon,p_{\pm},\Omega)\|g\|^{\frac{p_{-}}{p_{-}-1}}_{L^{q(\cdot)}(\Omega)},\ \text{or},\\
\varepsilon\|u\|^{p_{+}}_{L^{p(\cdot)}(\Omega)}+
C_{6}(\varepsilon,p_{\pm},\Omega)
\|g\|^{\frac{p_{+}}{p_{+}-1}}_{L^{q(\cdot)}(\Omega)},
\end{cases}\label{3.5}
\end{align}
where in the last inequality we used Young' inequality and
$\varepsilon\in(0,1)$ will be chosen later.

 Now we consider two
cases: (i) $\|\nabla u\|_{L^{p(\cdot)}(\Omega)}>1$, and (ii) $\|\nabla
u\|_{L^{p(\cdot)}(\Omega)}\leq 1$.

  (i) If $\|\nabla
u\|_{L^{p(\cdot)}(\Omega)}>1$, it follows from \eqref{2.2}, \eqref{3.2} and
\eqref{3.5} that
\begin{align}
\mathcal {J}_{\gamma}(u)\geq &  L^{-1}\int_{\Omega}|\nabla
u|^{p(x)}\text{d}x-\bigg|\int_{\Omega}gu\text{d}x\bigg|\label{3.6}\\
\geq &  ~L^{-1}\|\nabla
u\|^{p_{-}}_{L^{p(\cdot)}(\Omega)}-\bigg|\int_{\Omega}gu\text{d}x\bigg|\notag
\\
\geq &  ~L^{-1}C_{1}\|
u\|^{p_{-}}_{L^{p(\cdot)}(\Omega)}-L^{-1}\bigg(\|
\psi\|^{p_{-}}_{L^{p(\cdot)}(\Omega)}+\|\nabla\psi\|^{p_{-}}_{L^{p(\cdot)}(\Omega)}\bigg)-\varepsilon\|u\|^{p_{-}}_{L^{p(\cdot)}(\Omega)}-
C_{5}(\varepsilon,p_{\pm},\Omega)\|g\|^{\frac{p_{-}}{p_{-}-1}}_{L^{q(\cdot)}(\Omega)}.\label{3.7}
\end{align}
Choose $\varepsilon\in (0,1)$ such that $L^{-1}C_{1}-\varepsilon>0$,
then \eqref{3.7} yields
\begin{align*}\mathcal {J}_{\gamma}(u)>-L^{-1}\bigg(\|
\psi\|^{p_{-}}_{L^{p(\cdot)}(\Omega)}+\|\nabla\psi\|^{p_{-}}_{L^{p(\cdot)}(\Omega)}\bigg)-
C_{5}(\varepsilon,p_{\pm},\Omega)\|g\|^{\frac{p_{-}}{p_{-}-1}}_{L^{q(\cdot)}(\Omega)}>-\infty.
\end{align*}
(ii) If $\|\nabla u\|_{L^{p(\cdot)}(\Omega)}\leq 1$, we estimate by
\eqref{2.2}, \eqref{3.3} and \eqref{3.5}
\begin{align}
\mathcal {J}_{\gamma}(u)\geq &  ~L^{-1}\int_{\Omega}|\nabla
u|^{p(x)}\text{d}x-\bigg|\int_{\Omega}gu\text{d}x\bigg|\notag\\
\geq &  ~L^{-1}\|\nabla
u\|^{p_{+}}_{L^{p(\cdot)}(\Omega)}-\bigg|\int_{\Omega}gu\text{d}x\bigg|\notag
\\\geq &  ~L^{-1}C_{2}\|
u\|^{p_{+}}_{L^{p(\cdot)}(\Omega)}-L^{-1}\bigg(\|
\psi\|^{p_{+}}_{L^{p(\cdot)}(\Omega)}+\|\nabla\psi\|^{p_{+}}_{L^{p(\cdot)}(\Omega)}\bigg) -\varepsilon\|u\|^{p_{+}}_{L^{p(\cdot)}(\Omega)}-
C_{6}(\varepsilon,p_{\pm},\Omega)\|g\|^{\frac{p_{+}}{p_{+}-1}}_{L^{q(\cdot)}(\Omega)}.\label{3.8}
\end{align}
Choose $\varepsilon\in (0,1)$ such that $L^{-1}C_{2}-\varepsilon>0$,
then \eqref{3.8} gives
\begin{align*}\mathcal {J}_{\gamma}(u)>-L^{-1}\bigg(\|
\psi\|^{p_{+}}_{L^{p(\cdot)}(\Omega)}+\|\nabla\psi\|^{p_{+}}_{L^{p(\cdot)}(\Omega)}\bigg)-
C_{6}(\varepsilon,p_{\pm},\Omega)\|g\|^{\frac{p_{+}}{p_{+}-1}}_{L^{q(\cdot)}(\Omega)}>-\infty.
\end{align*}
 Let
us now prove existence of a minimizer of $\mathcal {J}_{\gamma}(u)$.
Let $u_{j}\in \mathcal {K}$ be a minimizing
sequence. We shall show that $\{u_{j}-\psi\}$ (up to a subsequence)
is bounded in $W^{1,p(\cdot)}_{0}(\Omega)$. Without loss of
generality, assume that $\|\nabla
u_{j}\|_{L^{p(\cdot)}(\Omega)}>1$ (If not, then $\|\nabla
u_{j}\|_{L^{p(\cdot)}(\Omega)}\leq 1$, which implies $\|
u_{j}-\psi\|_{L^{p(\cdot)}(\Omega)}\leq C\| \nabla
u_{j}-\nabla\psi\|_{L^{p(\cdot)}(\Omega)}\leq
C+C\|\nabla\psi\|_{L^{p(\cdot)}(\Omega)}<\infty$). Now for $j\gg 1$,
$\mathcal {J}_{\gamma}(u_{j}) \leq I_{0}+1$. From \eqref{3.6}, \eqref{3.4} and
\eqref{3.1} and applying Young' inequality with $\varepsilon$, we derive
\begin{align}
\|\nabla
u_{j}\|^{p_{-}}_{L^{p(\cdot)}(\Omega)}\leq &~\int_{\Omega}|\nabla
u_{j}|^{p(x)}\text{d}x\notag\\
\leq &~L\mathcal
{J}_{\gamma}(u_{j})+L\bigg|\int_{\Omega}gu_{j}\text{d}x\bigg|\notag\\
\leq &~L(I_{0}+1)+LC_{7}(p_{\pm},\Omega,\|g\|_{L^{q(\cdot)}(\Omega)})\|u_{j}\|_{L^{p(\cdot)}(\Omega)},\notag
\\\leq &~C_{8}(\|\nabla
u\|_{L^{p(\cdot)}(\Omega)}+\|\nabla\psi\|_{L^{p(\cdot)}(\Omega)}+
\|\psi\|_{L^{p(\cdot)}(\Omega)})+L(I_{0}+1),\notag\\
\leq &~\frac{1}{2}\|\nabla
u_{j}\|^{p_{-}}_{L^{p(\cdot)}(\Omega)}+C_{9}(1+\|\nabla\psi\|_{L^{p(\cdot)}(\Omega)}+
\|\psi\|_{L^{p(\cdot)}(\Omega)}),\notag
\end{align}
where $C_{8},C_{9}$ depend only on
$L,I_{0},p_{\pm},\Omega,\|g\|_{L^{q(\cdot)}(\Omega)}$. Therefore, we
get
\begin{align*}
|\nabla u_{j}\|^{p_{-}}_{L^{p(\cdot)}(\Omega)}\leq
2C_{9}(1+\|\nabla\psi\|_{L^{p(\cdot)}(\Omega)}+
\|\psi\|_{L^{p(\cdot)}(\Omega)}).
\end{align*}
Thus, using Poincar\'{e} inequaltiy once more, we deduce that
$\{u_{j}-\psi\}$ is bounded in $W^{1,p(\cdot)}_{0}(\Omega)$. By
reflexivity, there is a function $u\in
\mathcal {K}$ such that, up to a subsequence,
\begin{align}
u_{j}\rightharpoonup u \ \text{weakly\ in} \
W^{1,p(\cdot)}(\Omega),\ \ u_{j}\rightarrow u \ \text{in}\
L^{p(\cdot)}(\Omega),\ \ \ u_{j}\rightarrow u\ \text{a.e.\ in}\
\Omega. \notag
\end{align}
With a slight modification of \cite[Theorem 1.6]{LU}, we
deduce from \eqref{2.1} and \eqref{2.2} that
\begin{align}\label{3.10}
\int_{\Omega}f(x,|\nabla u|)\text{d}x\leq \liminf_{j\rightarrow
\infty}\int_{\Omega}f(x,|\nabla u_{j}|)\text{d}x.
\end{align}
By pointwise convergence we have, in the case of $0<\gamma\leq 1$,
\begin{align}\label{3.11}
\int_{\Omega}(F_{\gamma}(u)+gu)\text{d}x\leq \liminf_{j\rightarrow
\infty}\int_{\Omega}(F_{\gamma}(u_{j})+gu_{j})\text{d}x.
\end{align}

 For $\gamma=0$, recalling that $\lambda_{+}>\lambda_{-}>0$, we have
\begin{align}
\int_{\Omega}\lambda_{-}\chi_{\{u\leq 0\}}\text{d}x=&
\int_{\{u\leq 0\}}\lambda_{-}\chi_{\{u_{j}>
0\}}\text{d}x+\int_{\{u\leq
0\}}\lambda_{-}\chi_{\{u_{j}\leq 0\}}\text{d}x\notag\\
\leq &~\int_{\{u\leq 0\}}\lambda_{+}\chi_{\{u_{j}>
0\}}\text{d}x+\int_{\Omega}\lambda_{-}\chi_{\{u_{j}\leq
0\}}\text{d}x,\notag
\end{align}
which implies
\begin{align}
\int_{\Omega}\lambda_{-}\chi_{\{u\leq 0\}}\text{d}x
\leq\liminf_{j\rightarrow \infty}\bigg(\int_{\{u\leq
0\}}\lambda_{+}\chi_{\{u_{j}>
0\}}\text{d}x+\int_{\Omega}\lambda_{-}\chi_{\{u_{j}\leq
0\}}\text{d}x\bigg).\notag
\end{align}
On the other hand, since $u_{j}\rightarrow u$ a.e. in $\Omega$, it
follows from the Dominated Convergence Theorem that
\begin{align} \int_{\Omega}\lambda_{+}\vartheta_{\{u>
0\}}\text{d}x =& \int_{\{u> 0\}}\lambda_{+}(\lim_{j\rightarrow
\infty}\chi_{\{u_{j}> 0\}})\text{d}x\notag\\
=& \lim_{j\rightarrow \infty}\int_{\{u>
0\}}\lambda_{+}\chi_{\{u_{j}> 0\}}\text{d}x.\notag
\end{align}
Hence
\begin{align}\label{3.12}
\int_{\Omega}(F_{0}(u)+gu)\text{d}x\leq \liminf_{j\rightarrow
\infty}\int_{\Omega}(F_{0}(u_{j})+gu_{j})\text{d}x.
\end{align}

 Now from \eqref{3.10},\eqref{3.11} and \eqref{3.12} we conclude that
\begin{align}
\mathcal {J}_{\gamma}(u)\leq \liminf_{j\rightarrow \infty} \mathcal
{J}_{\gamma}(u_{j})=I_{0},\notag
\end{align}
for all $0\leq \gamma\leq 1$, which proves the existence of a
minimizer
under the condition of $g\in L^{q(\cdot)}(\Omega)$.

 Secondly, we establish the $L^{\infty}-$boundedness of $u_{\gamma}$,
provided $g\in L^{q(\cdot)}(\Omega)$. Hereafter in this proof we
will refer $u_{\gamma}$ as $u$.

 Let $j_{0}:=\big(\sup\limits_{\partial
\Omega} \psi\big)$ be the smallest natural number above $\sup\limits_{\partial
\Omega} \psi$. For each $j\geq j_{0}$, we define the truncated
function $u_{j}:\Omega\rightarrow\mathbb{R}$ by
\begin{align}
u_{j}=\left\{
  \begin{array}{ll}
    j\cdot sing(u),\ \ \text{if}\ |u|>j, \\
    u,\ \ \text{if}\ |u|\leq j,
  \end{array}
\right.\notag
\end{align}
where $sing(u)=1$ if $u\geq 0$ and $sing(u)=-1$ if $u< 0$. Define
the set $A_{j}:=\{|u|>j\}$. For $0<\gamma \leq 1$, in view of the
minimality of $u$, we derive
\begin{align}
\int_{A_{j}}f(x,\nabla u)\text{d}x=& \int_{\Omega}(f(x,\nabla u)
-f(x,\nabla u_{j}))+\int_{A_{j}}f(x,\nabla u_{j})\text{d}x\notag\\
\leq &~\int_{A_{j}}g(u_{j}-u)\text{d}x+\int_{A_{j}}
\lambda_{+}((u_{j}^{+})^{\gamma}-(u^{+})^{\gamma})\text{d}x+\int_{A_{j}}\lambda_{-}((u_{j}^{-})^{\gamma}-(u^{-})^{\gamma})\text{d}x+L|A_{j}|.\label{3.13}
\end{align}
 Now we estimate each integration in the right side of \eqref{3.13}.
\begin{align}
\int_{A_{j}}\lambda_{+}((u_{j}^{+})^{\gamma}-(u^{+})^{\gamma})\text{d}x
=& \lambda_{+}\int_{A_{j}\cap\{u>0\}}(j^{\gamma}-|u|^{\gamma})\text{d}x+\lambda_{+}\int_{A_{j}\cap\{u\leq
0\}}\bigg(((-j)^{+})^{\gamma}-(u^{+})^{\gamma}\bigg)\text{d}x\notag\\
\leq &~0.\notag
\end{align}
\begin{align}
\int_{A_{j}}\lambda_{-}((u_{j}^{-})^{\gamma}-(u^{-})^{\gamma})\text{d}x
=& \lambda_{-}\int_{A_{j}\cap\{u\leq 0\}}
(j^{\gamma}-|u|^{\gamma})\text{d}x+\lambda_{-}\int_{A_{j}\cap\{u>
0\}}\bigg((j^{-})^{\gamma}-(u^{-})^{\gamma}\bigg)\text{d}x\notag\\
\leq &~0.\notag
\end{align}
Then we find
\begin{align}
\int_{A_{j}}(F_{\gamma}(u_{j})-F_{\gamma}(u))\text{d}x\leq 0.\label{3.14}
\end{align}

 For the first integration in the right side of \eqref{3.13},
it follows
\begin{align}
\int_{A_{j}}g(u_{j}-u)\text{d}x=& \int_{A_{j}\cap\{u>0\}}g(j-u)\text{d}x+
\int_{A_{j}\cap\{u\leq 0\}}g(u-j)\text{d}x\notag\\
\leq &~2 \int_{A_{j}} |g|(|u|-j)\text{d}x.\label{3.15}
\end{align}
 For $\gamma =0$ it suffices to notice that $u_{j}>0$ and $u$ have
the same sign. By the choice of the truncated function, we know
that $(|u|-j)^{+}\in W^{1,p(\cdot)}_{0}(A_{j})$. Let
$
\frac{1}{t(\cdot)}=1-\frac{1}{p(\cdot)}-\frac{1}{q(\cdot)},\
t_{-}=\inf \limits_{x\in\Omega}t(x),\ t_{+}=\sup\limits _{x\in\Omega}t(x)$ and $p^{*}(\cdot)=\frac{np(\cdot)}{n-p(\cdot)}.\notag
$
Applying
H\"{o}lder's inequality and embedding theorem, we find
\begin{align}
\int_{A_{j}} |g|(|u|-j)^{+}\text{d}x\leq &~
2\|g\|_{L^{\frac{p(\cdot)}{p(\cdot)-1}}(A_{j})}
\|(|u|-j)^{+}\|_{L^{p(\cdot)}(A_{j})}\notag\\
\leq &~C\|g\|_{L^{q(\cdot)}(A_{j})}\|1\|_{L^{t(\cdot)}(A_{j})}
\|(|u|-j)^{+}\|_{L^{p^{*}(\cdot)}(A_{j})}\|1\|_{L^{n}(A_{j})}\notag\\
\leq &~\left\{
  \begin{array}{ll}
    C\|g\|_{L^{q(\cdot)}(\Omega)}|A_{j}|^{
\frac{1}{t_{-}}+\frac{1}{n}} \|\nabla(|u|-j)^{+}\|_{L^{p(\cdot)}(A_{j})},\ \ \text{if}
\ \ |A_{j}|>1\\
    C\|g\|_{L^{q(\cdot)}(\Omega)}|A_{j}|^{
\frac{1}{t_{+}}+\frac{1}{n}}
\|\nabla(|u|-j)^{+}\|_{L^{p(\cdot)}(A_{j})},\ \ \text{if} \ \
|A_{j}|\leq 1
  \end{array}
\right.\notag\\
=& \left\{
  \begin{array}{ll}
    C|\Omega|^{
\frac{1}{t_{-}}+\frac{1}{n}}(\frac{|A_{j}|}{|\Omega|})^{
\frac{1}{t_{-}}+\frac{1}{n}}
\|\nabla(|u|-j)^{+}\|_{L^{p(\cdot)}(A_{j})},\ \ \text{if}
\ \ |A_{j}|>1\\
    C|\Omega|^{
\frac{1}{t_{+}}+\frac{1}{n}}(\frac{|A_{j}|}{|\Omega|})^{
\frac{1}{t_{+}}+\frac{1}{n}}
\|\nabla(|u|-j)^{+}\|_{L^{p(\cdot)}(A_{j})},\ \ \text{if} \ \
|A_{j}|\leq 1
  \end{array}
\right.\notag\\
 \leq &~C\bigg(1+|\Omega|^{
\frac{1}{t_{-}}+\frac{1}{n}}\bigg)\left(\frac{|A_{j}|}{|\Omega|}\right)^{
\frac{1}{t_{+}}+\frac{1}{n}}\|\nabla u\|_{L^{p(\cdot)}(A_{j})}\notag\\
=& C\left(\frac{|A_{j}|}{|\Omega|}\right)^{\frac{1}{t_{+}}+\frac{1}{n}}\|\nabla
u \|_{L^{p(\cdot)}(A_{j})},\label{3.16}
\end{align}
where the constant $C$ in the last inequality depends only on
$p_{\pm},q_{-},n,\Omega,\|g\|_{L^{q(\cdot)}(\Omega)}$.

 Collecting \eqref{3.13}-\eqref{3.16}, we obtain
\begin{align}
\int_{A_{j}}f(x,\nabla u)\text{d}x\leq
C\left(\frac{|A_{j}|}{|\Omega|}\right)^{\frac{1}{t_{+}}+\frac{1}{n}}\|\nabla
u \|_{L^{p(\cdot)}(A_{j})}+L|A_{j}|,\label{3.17}
\end{align}
where $C$ depends only on
$p_{\pm},q_{-},n,\Omega,\|g\|_{L^{q(\cdot)}(\Omega)}$.

 Now we consider two cases: (i)
$\|\nabla u\|_{L^{p(\cdot)}(A_{j})}>1$, and (ii) $\|\nabla
u\|_{L^{p(\cdot)}(A_{j})}\leq 1$.

 (i) If $\|\nabla
u\|_{L^{p(\cdot)}(A_{j})}>1$, we estimate by \eqref{2.2}, \eqref{3.17} and Young' inequality
\begin{align}
\|\nabla u\|_{L^{p(\cdot)}(A_{j})}^{p_{-}}\leq &~\int_{A_{j}}|\nabla
u|^{p(x)}\text{d}x\notag\\
\leq &~L\int_{A_{j}}f(x,\nabla u)\text{d}x\notag\\
\leq &~C\left(\frac{|A_{j}|}{|\Omega|}\right)^{\frac{1}{t_{+}}+\frac{1}{n}}\|\nabla
u
\|_{L^{p(\cdot)}(A_{j})}+L^{2}|A_{j}|\notag\\
\leq &~C\left(\frac{|A_{j}|}{|\Omega|}\right)^{(\frac{1}{t_{+}}+\frac{1}{n})\frac{p_{-}}{p_{-}-1}}+
\frac{1}{2}\|\nabla
u\|_{L^{p(\cdot)}(A_{j})}^{p_{-}}+L^{2}|A_{j}|,\notag
\end{align}
which implies
\begin{align}
\|\nabla u\|_{L^{p(\cdot)}(A_{j})}^{p_{-}}\leq
C\left(\frac{|A_{j}|}{|\Omega|}\right)^{(\frac{1}{t_{+}}+\frac{1}{n})\frac{p_{-}}{p_{-}-1}}+L^{2}|A_{j}|
=C\left(\frac{|A_{j}|}{|\Omega|}\right)^{(1-\frac{1}{p_{-}}-\frac{1}{q_{-}}+\frac{1}{n})\frac{p_{-}}{p_{-}-1}}+L^{2}|A_{j}|.\notag
\end{align}
Therefore
\begin{align}
\|\nabla u\|_{L^{p(\cdot)}(A_{j})}\leq
C\left(\frac{|A_{j}|}{|\Omega|}\right)^{(1-\frac{1}{p_{-}}-\frac{1}{q_{-}}+\frac{1}{n})
\frac{1}{p_{-}-1}}+C\left(\frac{|A_{j}|}{|\Omega|}\right)^{\frac{1}{p_{-}}},\label{3.18}
\end{align}
 where $C$ depends only on
$L,p_{\pm},q_{-},n,\Omega,\|g\|_{L^{q(\cdot)}(\Omega)}$.

 On the
other hand, by an analogue argument as \eqref{3.16} and Young' inequality,
we obtain
\begin{align}
\int_{A_{j}}(|u|-j)^{+}\text{d}x\leq &~
2\|1\|_{L^{\frac{p(\cdot)}{p(\cdot)-1}}(A_{j})}
\|(|u|-j)^{+}\|_{L^{p(\cdot)}(A_{j})}\notag\\
\leq &~\left\{
  \begin{array}{ll}
    C|A_{j}|^{1-
\frac{1}{p_{+}}+\frac{1}{n}} \|\nabla u\|_{L^{p(\cdot)}(A_{j})},\ \
\text{if}
\ \ |A_{j}|>1\\
    C|A_{j}|^{1-
\frac{1}{p_{-}}+\frac{1}{n}} \|\nabla u\|_{L^{p(\cdot)}(A_{j})},\ \
\text{if} \ \ |A_{j}|\leq 1
  \end{array}
\right.\notag\\
 \leq &~C\left(\frac{|A_{j}|}{|\Omega|}\right)^{1-
\frac{1}{p_{-}}+\frac{1}{n}}\|\nabla u\|_{L^{p(\cdot)}(A_{j})}\notag\\
\leq &~C\left(\frac{|A_{j}|}{|\Omega|}\right)^{1-
\frac{1}{p_{-}}+\frac{1}{n}}\bigg(\left(\frac{|A_{j}|}{|\Omega|}\right)^{(1-\frac{1}{p_{-}}-\frac{1}{q_{-}}+\frac{1}{n})
\frac{1}{p_{-}-1}}+C\left(\frac{|A_{j}|}{|\Omega|}\right)^{\frac{1}{p_{-}}}\bigg)\notag\\
=&  C\left(\frac{|A_{j}|}{|\Omega|}\right)^{(1-
\frac{1}{p_{-}}-\frac{1}{q_{-}}+\frac{1}{n})\frac{1}{p_{-}-1}+(1-
\frac{1}{p_{-}}-\frac{1}{n})}+C\left(\frac{|A_{j}|}{|\Omega|}\right)^{1+\frac{1}{n}},\label{3.19}
\end{align}
where in the last inequality we used \eqref{3.18}, the constant $C$ depends
only on $L,p_{\pm},q_{-},n,\Omega,\|g\|_{L^{q(\cdot)}(\Omega)}$.

 (ii)
If $\|\nabla u\|_{L^{p(\cdot)}(A_{j})}\leq 1$, analogously, we
deduce that
\begin{align}
\int_{A_{j}}(|u|-j)^{+}\text{d}x\leq
C\left(\frac{|A_{j}|}{|\Omega|}\right)^{(1-
\frac{1}{p_{-}}-\frac{1}{q_{-}}+\frac{1}{n})\frac{1}{p_{+}-1}+(1-
\frac{1}{p_{-}}-\frac{1}{n})}+C\left(\frac{|A_{j}|}{|\Omega|}\right)^{1+\frac{1}{n}},\label{3.20}
\end{align}
where the constant $C$ depends only on
$L,p_{\pm},q_{-},n,\Omega,\|g\|_{L^{q(\cdot)}(\Omega)}$.\\
Now combining \eqref{3.19} and \eqref{3.20}, we get
\begin{align*}
\int_{A_{j}}(|u|-j)^{+}\text{d}x \leq
C\left(\frac{|A_{j}|}{|\Omega|}\right)^{(1-
\frac{1}{p_{-}}-\frac{1}{q_{-}}+\frac{1}{n})\frac{1}{p_{+}-1}+(1-
\frac{1}{p_{-}}-\frac{1}{n})}+C\left(\frac{|A_{j}|}{|\Omega|}\right)^{1+\frac{1}{n}},
\end{align*}
where $\epsilon_{0}=\min\{\frac{1}{n},(1-
\frac{1}{p_{-}}-\frac{1}{q_{-}}+\frac{1}{n})\frac{1}{p_{+}-1}+(1-
\frac{1}{p_{-}}-\frac{1}{n})-1\}$ and $C$ depends only on
$L,p_{\pm},q_{-},n,\Omega,\|g\|_{L^{q(\cdot)}(\Omega)}$. Notice that
by \eqref{2.7} we have
$\frac{1}{q_{-}}<\frac{1}{n}-\frac{1}{p_{-}}+\frac{1}{p_{+}}$, thus
$\epsilon_{0}>0$. Notice also that $\|u\|_{L^{1}(A_{j_{0}})}\leq
\bigg(1+|A_{j_{0}}|^{\frac{p_{-}-1}{p_{-}}}\bigg)\|u\|_{L^{p(x)}(A_{j_{0}})}\leq
C$.
Applying \cite[Lemma 5.1]{S}, we obtain the desired result. \hfill $\blacksquare$
\end{pf*}
\begin{remark}\label{Remark 3.2} Note that in \cite{EH3}, the assumption that $\int_{\Omega} |\nabla u|^{p(x)}\text{d}x\leq M$ with some $M\geq 0$ is assumed in the establishment of local regularity for minimizers of a functional with a form $\int_{\Omega}f(x,u,\nabla u)\text{d}x$, while in this paper, we can show that any minimizer $u_{\gamma}$ of $\mathcal {J}_{\gamma}(u)$ governed by \eqref{1.1} is uniformly bounded in $W^{1,p(\cdot)}(\Omega)$ by $L^{\infty}-$estimates of $u_{\gamma}$. Indeed,
we have
\begin{align}
\int_{\Omega} |\nabla u_\gamma|^{p(x)}\text{d}x\leq &~ L\int_{\Omega}
f(x,\nabla u_\gamma)\text{d}x\notag\\
\leq &~ L\bigg(\mathcal{J}_{\gamma}(\psi)-\int_{\Omega} F(
u_\gamma)\text{d}x+\int_{\Omega} |gu_\gamma|\text{d}x\bigg)\notag\\
\leq &  L\mathcal{J}_{\gamma}(\psi)+C(L,n,p_{\pm},
\lambda_{\pm},\Omega,\|\psi\|_{L^{\infty}(\partial\Omega)},
\|g\|_{L^{q(\cdot)}(\Omega)})\notag\\
\leq &~ M,\notag
\end{align}
where $M=M(L,n,q_{-},p_{\pm},
\lambda_{\pm},\Omega,\|\psi\|_{L^{\infty}(\partial\Omega)},
\|g\|_{L^{q(\cdot)}(\Omega)})$ is a positive constant. Therefore, we
conclude by $ u_{\gamma}-\psi\in W^{1,p(\cdot)}_0(\Omega)$ that
$
\|u_{\gamma}\|_{W^{1,p(\cdot)}(\Omega)} \leq  C,\notag
$ where $C $ is independent of $\gamma$.
\end{remark}
\section{High integrability}
In this section we prove a higher integrability result for minimizers of functional in \eqref{1.1}.
\begin{proposition}\label{Proposition 4.1}
Assume that \eqref{2.1}-\eqref{2.7} hold. Let
 $u\in \mathcal {K}$ be a
minimizer of the functional $\mathcal {J}_{\gamma}$ governed by
\eqref{1.1}. Then there exist two positive constants $C_{0}$ and
$\delta_{0}<q_{-}(1-\frac{1}{p_{-}})-1$, both depending only on
$n,p_{\pm},\lambda_{\pm},q_{-},L,M,\Omega$, such that
\begin{align}
\bigg(\frac{1}{|B_{R/2}|}\int_{B_{R/2}}|\nabla
u|^{p(x)(1+\delta_{0})}\text{d}x\bigg)^{\frac{1}{1+\delta_{0}}}\leq &~ \frac{C_{0}}{|B_{R}|}\int_{B_{R}}|\nabla
u|^{p(x)}\text{d}x+C_{0}
\bigg(\frac{1}{|B_{R}|}\int_{B_{R}}(1+|g|^{\frac{p_{-}}{p_{-}-1}(1+\delta_{0})})\text{d}x\bigg)
^{\frac{1}{1+\delta_{0}}},\label{4.1}
\end{align}
for all $B_{R}\Subset\Omega$.
\end{proposition}
  In order to prove Proposition 4.1, we need the following iteration lemma.
  \begin{lemma}\label{Lemma 4.2}\cite{EH3}
Let $0<\theta<1,\ A>0,\ B\geq 0,\
1<p_{-}\leq p(x)\leq p_{+}<+\infty$, and let $f\geq 0$ be a bounded
function on $(r,R)$ satisfying
\begin{align}
f(t)\leq \theta
f(s)+A\int_{B_{R}}\left|\frac{h(x)}{s-t}\right|^{p(x)}\text{d}x+B,\notag
\end{align}
for all $r\leq t<s\leq R$, where $h\in L^{p(\cdot)}(B_{R})$. Then
there exists a constant $C= C(\theta, p_{+})$ such that
\begin{align}
f(r)\leq
C\bigg(A\int_{B_{R}}\left|\frac{h(x)}{R-r}\right|^{p(x)}\text{d}x+B\bigg).\notag
\end{align}
\end{lemma}
\begin{pf*}{Proof of Proposition~\ref{Proposition 4.1}}
Let $0<R<R_{0}\leq 1$ and let
$x_{0}\in B_{R}$ with $\overline{B}_{R_{0}}(x_{0})\subset\Omega$.
Let $t,s\in\mathbb{R}$ with $\frac{R}{2}<t<s<R$. Let $\eta\in
C_{c}^{\infty}(B_{R}),0\leq \eta\leq 1$, be a cut-off function with
$\eta\equiv1$ on $B_{t},\eta\equiv0$ outside $B_{s}$ and $|\nabla
\eta|\leq \frac{2}{s-t}$. We define the function
$z=u-\eta(u-(u)_{R})$. We deduce from \eqref{2.2} and minimality of $u$
that
\begin{align}
L^{-1}\int_{B_{t}}|\nabla u|^{p(x)}\text{d}x\leq &~
\int_{B_{t}}f(x,\nabla u)\text{d}x\notag\\
\leq &~
\int_{B_{s}}f(x,\nabla u)\text{d}x\notag\\
\leq &~\int_{B_{s}}f(x,\nabla
z)+(F_{\gamma}(z)-F_{\gamma}(u))+g(z-u)\text{d}x\notag\\
\leq &~L\int_{B_{s}}(\mu^{2}+|\nabla
z|^{2})^{\frac{p(x)}{2}}\text{d}x+\int_{B_{s}}F_{\gamma}(z)
-F_{\gamma}(u)\text{d}x+\int_{B_{s}} g(z-u)\text{d}x,\label{4.2}
\end{align}
where in the last but one inequality we used the fact that if
$\varphi\in W^{1,p(\cdot)}_{0}(\Omega)$ with spt
$\varphi\Subset\Omega$, then there holds
\begin{align}
\int_{\text{spt}\ \varphi}\big(f(x,\nabla
u)+F_{\gamma}(u)+gu\big)\text{d}x\leq \int_{\text{spt}\
\varphi}\big(f(x,\nabla u+\nabla
\varphi)+F_{\gamma}(u+\varphi)+g(u+\varphi)\big)\text{d}x.\notag
\end{align}
Indeed, it follows from the minimality of $u$ that
\begin{align}
\int_{\text{spt} \varphi}&\big(f(x,\nabla
u)+F_{\gamma}(u)+gu\big)\text{d}x+\int_{\Omega\setminus(\text{spt}\
\varphi)}\big(f(x,\nabla u)+F_{\gamma}(u)+gu\big)\text{d}x\notag\\
 \leq &~ \int_{\text{spt}\ \varphi}\big(f(x,\nabla u+\nabla
\varphi)+F_{\gamma}(u+\varphi)+g(u+\varphi)\big)\text{d}x+\int_{\Omega\setminus(\text{spt}\ \varphi)}\big(f(x,\nabla
u+\nabla
\varphi)+F_{\gamma}(u+\varphi)+g(u+\varphi)\big)\text{d}x\notag\\
\leq &~\int_{\text{spt}\ \varphi}\big(f(x,\nabla u+\nabla
\varphi)+F_{\gamma}(u+\varphi)+g(u+\varphi)\big)\text{d}x+\int_{\Omega\setminus(\text{spt}\ \varphi)}\big(f(x,\nabla
u+\nabla
\varphi)+F_{\gamma}(u+\varphi)+g(u+\varphi)\big)\text{d}x\notag\\
=& \int_{\text{spt}\ \varphi}\big(f(x,\nabla u+\nabla
\varphi)+F_{\gamma}(u+\varphi)+g(u+\varphi)\big)\text{d}x+
\int_{\Omega\setminus(\text{spt}\ \varphi)}\big(f(x,\nabla
u)+F_{\gamma}(u)+gu\big)\text{d}x.\notag
\end{align}
 We shall estimate each integration of \eqref{4.2}.
\begin{align}
\int_{B_{s}}|\nabla z|^{p(x)}\text{d}x\leq &~
\int_{B_{s}}|(1-\eta)\nabla u-\nabla \eta
(u-(u)_{R})|^{p(x)}\text{d}x\notag\\
\leq &~C\int_{B_{s}\setminus B_{t}}|\nabla
u|^{p(x)}\text{d}x+C\int_{B_{s}}\left|\frac{u-(u)_{R}}{s-t}\right|^{p(x)}\text{d}x,\label{4.3}
\end{align}
where $C=C(p_{+},p_{-})$ is a positive constant.

 A direct calculus shows that
\begin{align}
\int_{B_{s}}F_{\gamma}(z) -F_{\gamma}(u)\text{d}x=& \lambda_{+}\int_{B_{s}}
((z^{+})^{\gamma}-(u^{+})^{\gamma})\text{d}x
+\lambda_{-}\int_{B_{s}}
((z^{-})^{\gamma}-(u^{-})^{\gamma})\text{d}x\notag\\
\leq &~C\int_{B_{s}}|z-u|^{\gamma}\text{d}x,\notag
\end{align}
where $C=C(\lambda_{+},\lambda_{-})$ is a positive constant.

 Then we estimate from Young' inequality that
\begin{align}
\int_{B_{s}}F_{\gamma}(z) -F_{\gamma}(u)\text{d}x \leq &~C\int_{B_{s}}|u-(u)_{R}|^{\gamma}\text{d}x=
C\int_{B_{s}}\left|\frac{u-(u)_{R}}{s-t}\right|^{\gamma}|s-t|^{\gamma}\text{d}x\notag\\
\leq &~C\int_{B_{s}}\left|\frac{u-(u)_{R}}{s-t}\right|^{p(x)}\text{d}x+C\int_{B_{s}}|s-t|^{\frac{\gamma
p(x)}{p(x)-\gamma}}\text{d}x\notag\\
=& C\int_{B_{s}}\left|\frac{u-(u)_{R}}{s-t}\right|^{p(x)}\text{d}x+C|B_{s}|,\label{4.4}
\end{align}
where $C=C(p_{\pm},\lambda_{\pm})$ is a positive constant.
\begin{align}
\int_{B_{s}}|g(z-u)|\text{d}x\leq &~\int_{B_{s}}|g||u-(u)_{R}|\text{d}x\notag\\
\leq &~C\int_{B_{s}}\left|\frac{u-(u)_{R}}{s-t}\right|^{p(x)}\text{d}x+
C\int_{B_{s}}(|g||s-t|)^{\frac{p(x)}{p(x)-1}}\text{d}x\notag\\
\leq &~C\int_{B_{s}}\left|\frac{u-(u)_{R}}{s-t}\right|^{p(x)}\text{d}x+
C\int_{B_{s}}|g|^{\frac{p(x)}{p(x)-1}}\text{d}x\notag\\
\leq &~C\int_{B_{s}}\left|\frac{u-(u)_{R}}{s-t}\right|^{p(x)}\text{d}x+
C\int_{B_{s}}\bigg(1+|g|^{\frac{p_{-}}{p_{-}-1}}\bigg)\text{d}x,\label{4.5}
\end{align}
where $C=C(p_{+},p_{-})$ is a positive constant.

 Combining \eqref{4.2}-\eqref{4.5}, we obtain
\begin{align*}
\int_{B_{t}}|\nabla u|^{p(x)}\text{d}x\leq
C\int_{B_{s}\setminus B_{t}}|\nabla
u|^{p(x)}\text{d}x+C\int_{B_{s}}\left|\frac{u-(u)_{R}}{s-t}\right|^{p(x)}\text{d}x+
C\int_{B_{s}}\bigg(1+|g|^{\frac{p_{-}}{p_{-}-1}}\bigg)\text{d}x,
\end{align*}
where the constant $C$ depends only on $L,p_{\pm},\lambda_{\pm}$.

 Now ``filling the hole", we get
\begin{align*}
\int_{B_{t}}|\nabla u|^{p(x)}\text{d}x\leq
\frac{C}{1+C}\int_{B_{s}}|\nabla
u|^{p(x)}\text{d}x+\int_{B_{s}}\left|\frac{u-(u)_{R}}{s-t}\right|^{p(x)}\text{d}x+
\int_{B_{s}}\bigg(1+|g|^{\frac{p_{-}}{p_{-}-1}}\bigg)\text{d}x,\notag
\end{align*}
which and Lemma 4.2 imply
\begin{align}
\frac{1}{|B_{R/2}|}\int_{B_{R/2}}|\nabla
u|^{p(x)}\text{d}x\leq
C\frac{1}{|B_{R}|}\int_{B_{R}}\left|\frac{u-(u)_{R}}{R-R/2}\right|^{p(x)}\text{d}x+
C\frac{1}{|B_{R}|}\int_{B_{R}}\bigg(1+|g|^{\frac{p_{-}}{p_{-}-1}}\bigg)\text{d}x.\label{4.7}
\end{align}
Let
$
p_{1}=\min\limits_{x\in \overline{B}_{R}}p(x),\ \  p_{2}=\max\limits_{x\in
\overline{B}_{R}}p(x). \notag
$
By Sobolev-Poincar\'{e}'s
inequality, there exists $\nu<1$ such that
\begin{align}
\frac{1}{|B_{R}|}\int_{B_{R}}\left|\frac{u-(u)_{R}}{R}\right|^{p(x)}\text{d}x\leq &~1+\frac{1}{|B_{R}|}\int_{B_{R}}\left|\frac{u-(u)_{R}}{R}\right|^{p_{2}}\text{d}x\notag\\
\leq &~1+C\bigg(\int_{B_{R}}(1+|\nabla
u|^{p(x)})\text{d}x\bigg)^{\frac{p_{2}-p_{1}}{p_{1}\nu}}
R^{\frac{(p_{1}-p_{2})n}{p_{1}\nu}}\bigg(\frac{1}{|B_{R}|}\int_{B_{R}}|\nabla
u|^{p_{1}\nu}\text{d}x\bigg)^{\frac{1}{\nu}}\notag\\
\leq &~
C\left(\frac{1}{|B_{R}|}\int_{B_{R}}|\nabla
u|^{p(x)\nu}\text{d}x\right)^{\frac{1}{\nu}}+C,\label{4.8}
\end{align}
where in the last inequality we used Remark \ref{Remark 3.2} and the fact that, by
\eqref{2.5}, $R^{\frac{(p_{1}-p_{2})n}{p_{1}\nu}}$ is bounded.\\
Combining \eqref{4.7} and \eqref{4.8}, we get
\begin{align}
\frac{1}{|B_{R/2}|}\int_{B_{R/2}}|\nabla
u|^{p(x)}\text{d}x\leq
C\bigg(\frac{1}{|B_{R}|}\int_{B_{R}}|\nabla
u|^{p(x)\nu}\text{d}x\bigg)^{\frac{1}{\nu}}+
C\frac{1}{|B_{R}|}\int_{B_{R}}\bigg(1+|g|^{\frac{p_{-}}{p_{-}-1}}\bigg)\text{d}x,\notag
\end{align}
where $C=C(n,p_{\pm},\lambda_{\pm},L,M,\Omega)$. We now apply Gehring's lemma (see \cite{DGK}) to deduce that there exists
$0<\delta_{0}<q_{1}(1-\frac{1}{p_{-}})-1$ such that \eqref{4.1} holds.\hfill $\blacksquare$
\end{pf*}
\section{H\"{o}lder estimates for minimizers of functional $\mathcal {H}_{\gamma}$}
In this section, we establish local $C^{0,\alpha}-$regularity for
minimizers of the functional $\mathcal
{H}_{\gamma}(\gamma\in[0,1])$ governed
by \eqref{2.9}. We always let $v\in W^{1,p(\cdot)}(B_{R})$ with $v-u\in W^{1,p(\cdot)}_{0}(B_{R})$ be a
minimizer of the following local integral functional
 \begin{align}
\mathcal {H}_{\gamma}(v)=\int_{B_{R}(x_{0})}\big(h(\nabla
v)+F_{\gamma}(v)+gv\big)\text{d}x,\ \ B_{R}(x_{0})\Subset \Omega,\label{5.1}
\end{align}
and let $\tilde{v}(y)=\frac{1}{R}v(x_{0}+Ry)$. It is easy to check that
$\tilde{v}$ is a minimizer of the functional
\begin{align}
\tilde{\mathcal {H}}_{\gamma}(\tilde{v})=\int_{B_{1}(0)}[h(\nabla
\tilde{v})+R^{\gamma}F_{\gamma}(\tilde{v})+Rg\tilde{v}]\text{d}y,\label{5.2}
\end{align}
in the class $\big\{ \tilde{v} \in W^{1,\tilde{p}(\cdot)}(B_{1}): \tilde{v}-\frac{u}{R}\in W^{1,\tilde{p}(\cdot)}_{0}(B_{1})\big\}$, where
$\tilde{p}(y)=p(x_{0}+Ry)$.

Let
$
p_{1}=\min\limits_{x\in \overline{B}_{R}(x_{0})}p(x),\ \  p_{2}=\max\limits_{x\in
\overline{B}_{R}(x_{0})}p(x). \notag
$

 The following lemma is a slight version of \cite[Lemma 7.1]{G}, and can be obtained by induction in the same way as in \cite[Lemma 7.1]{G}. We omit the proof here.
\begin{lemma}\label{Lemma 5.1}
Let $0<a_{1}\leq a_{2}$ and
$\{\vartheta_{i}\}$ be a sequence of real positive numbers, such
that
\begin{align}
\vartheta_{i+1}\leq
CB^{i}(\vartheta_{i}^{1+a_{1}}+\vartheta_{i}^{1+a_{2}}),\notag
\end{align}
with $C>1$ and $B>1$. If $\vartheta_{0}\leq
(2C)^{-\frac{1}{a_{1}}}B^{-\frac{1}{a_{1}^{2}}}$, then we have
\begin{align}
\vartheta_{i}\leq B^{-\frac{i}{a_{1}}}\vartheta_{0},\notag
\end{align}
and hence in particular $\lim\limits_{i\rightarrow\infty}\vartheta_{i}=0$.
\end{lemma}
%
\begin{lemma}\label{Lemma 5.2}
\cite{LdT} Let $\phi(s)$ be a non-negative and
non-decreasing function. Suppose that
\begin{align}
\phi(r)\leq
C_{1}\bigg(\left(\frac{r}{R}\right)^{\alpha}+\mu\bigg)\phi(R)+C_{2}R^{\beta},\notag
\end{align}
for all $r\leq R\leq R_{0}$, with $0<\beta<\alpha$, $ C_{1}$
positive constants and $C_{2},\mu$ non-negative constants. Then, for
any $\sigma\leq \beta$, there exists a constant
$\mu_{0}=\mu_{0}(C_{1},\alpha,\beta,\sigma)$ such that if
$\mu<\mu_{0}$, then for all $r\leq R\leq R_{0}$ it follows that
\begin{align}
\phi(r)\leq C_{3}r^{\sigma},\notag
\end{align}
where $C_{3}=C_{3}(C_{1},C_{2},R_{0},\phi,\sigma,\beta)$ is a
positive constant.
\end{lemma}
\begin{lemma}\label{Lemma 5.3}
If $\tilde{v}$
 is a minimizer of $\tilde{\mathcal {H}}_{\gamma}$ governed by
 \eqref{5.2}, then $\tilde{v}$ is locally bounded and satisfies the
 estimates
\begin{align}
\sup\limits_{B_{\frac{1}{2}}(0)}|\tilde{v}|\leq
C\bigg(\left(\int_{B_{1}(0)}|\tilde{v}|^{p_{2}}\text{d}y\right)^{\frac{1}{p_{2}}}+1\bigg),\label{5.3}
\end{align}
and
\begin{align}
\sup\limits_{B_{\frac{1}{2}}(0)}\tilde{v}\leq
C\bigg(\left(\int_{B_{1}(0)}(\tilde{v}^{+})^{p_{2}}\text{d}y\right)^{\frac{1}{p_{2}}}|A_{0,1}|^{\frac{\alpha}{p_{2}}}+1\bigg),\label{5.4}
\end{align}
for some $\alpha>0$, where $A_{0,1}=\{y\in
B_{1}(0):\tilde{v}(y)>0\}$,
$C=C(n,L,p_{\pm},\lambda_{\pm},q_{-},M,\Omega,\|g\|_{L^{q(\cdot)}(\Omega)})$
 is a positive constant.
 \end{lemma}
%
 \begin{pf*}{Proof of Lemma \ref{Lemma 5.3}}
Without loss of generality, we may
assume that $R\leq 1$. The proof proceeds in three steps.

 \emph{First step: De Giorgi type estimates}.\ \ For any
$k\in\mathbb{R}$, we define the sets
\begin{align}
A_{k,\sigma}=\{y\in B_{\sigma}(0):u(y)>k\},\ B_{k,\sigma}=\{y\in
B_{\sigma}(0):u(y)<k\}.\notag
\end{align}
 We claim that for any $k\in\mathbb{R},\tilde{v}$ satisfies the
inequalities
\begin{align}
\int_{A_{k,\sigma}}|\nabla
\tilde{v}(y)|^{\tilde{p}(y)}\text{d}y\leq &~
C_{1}\int_{A_{k,\tau}}\left|\frac{\tilde{v}(y)-k}{\tau-\sigma}\right|^{\tilde{p}(y)}\text{d}y
+C_{2}\int_{A_{k,\tau}}|\tau-\sigma|^{\frac{\gamma
p_{2}}{p_{2}-\gamma}}\text{d}y\notag\\
&+C_{2}\int_{A_{k,\tau}}\bigg(1+(|g||\tau-\sigma|)^{\frac{p_{1}}{p_{1}-1}}\bigg)\text{d}y,\label{5.5}
\end{align}
and
\begin{align}
\int_{B_{k,\sigma}}|\nabla
\tilde{v}(y)|^{\tilde{p}(y)}\text{d}y\leq &~
C_{1}\int_{B_{k,\tau}}\left|\frac{\tilde{v}(y)-k}{\tau-\sigma}\right|^{\tilde{p}(y)}\text{d}y
+C_{2}\int_{B_{k,\tau}}|\tau-\sigma|^{\frac{\gamma
p_{2}}{p_{2}-\gamma}}\text{d}y\notag\\
&+C_{2}\int_{B_{k,\tau}}\bigg(1+(|g||\tau-\sigma|)^{\frac{p_{1}}{p_{1}-1}}\bigg)\text{d}y,\label{5.6}
\end{align}
for any $\frac{1}{2}\leq \sigma<\tau\leq 1$, where
$C_{i}=C_{i}(L,\lambda_{\pm},p_{\pm})$. Indeed, for $\frac{1}{2}\leq
\sigma\leq s<t\leq\tau\leq 1$, let $\eta\in
C_{c}^{\infty}(B_{1}(0))$ with $\text{spt}\eta\subset
B_{t},\eta\equiv 1$ on $B_{s}(0),|\nabla \eta|\leq \frac{2}{t-s}$ be
a standard cut-off function. Set $\tilde{z}(y)=\tilde{v}(y)-\eta
\tilde{w}(y)$, where $\tilde{w}(y)=\max\{\tilde{v}(y)-k,0\}$. In
view of minimality of $\tilde{v}$, we obtain
\begin{align}
\int_{A_{k,s}}|\nabla
\tilde{v}(y)|^{\tilde{p}(y)}\text{d}y\leq &~
\int_{A_{k,t}}|\nabla \tilde{v}(y)|^{\tilde{p}(y)}\text{d}y\notag\\
\leq &~ C\bigg(\int_{A_{k,t}}|\nabla
\tilde{z}(y)|^{\tilde{p}(y)}\text{d}y+\int_{A_{k,t}}
R^{\gamma}(F_{\gamma}(\tilde{z})-F_{\gamma}(\tilde{v}))
+Rg(\tilde{z}-\tilde{v})\text{d}y\bigg)\notag\\
\leq &~C\bigg(\int_{A_{k,t}}|(1-\eta)\nabla
\tilde{v}-\nabla
\eta\cdot(\tilde{v}-k)|^{\tilde{p}(y)}\text{d}y+\int_{A_{k,t}}
(F_{\gamma}(\tilde{z})-F_{\gamma}(\tilde{v}))
+g(\tilde{z}-\tilde{v})\text{d}y\bigg)\notag\\
\leq &~C\int_{A_{k,t}\setminus A_{k,s}}|\nabla
\tilde{v}(y)|^{\tilde{p}(y)}\text{d}y+C\int_{A_{k,t}}
\left|\frac{\tilde{v}-k}{t-s}\right|^{\tilde{p}(y)}\text{d}y\notag\\
&+C\int_{A_{k,t}}
(F_{\gamma}(\tilde{z})-F_{\gamma}(\tilde{v}))
+g(\tilde{z}-\tilde{v})\text{d}y,\label{5.7}
\end{align}
where $C=C(\tilde{p}_{1},\tilde{p}_{2})$ is a positive constant. We
remark that $\tilde{p}_{1}=\min\limits_{y\in
\overline{B}_{1}(0)}\tilde{p}(y)=p_{1},\ \tilde{p}_{2}=\max\limits_{y\in
\overline{B}_{1}(0)}\tilde{p}(y)=p_{2}$. Therefore
$C=C(\tilde{p}_{1},\tilde{p}_{2})=C(p_{1},p_{2})$. Moreover, we can
let $C$ depend only on $p_{\pm}$.

 In view of \eqref{4.4} and \eqref{4.5}, we derive
\begin{align*}
&\int_{A_{k,t}}F_{\gamma}(\tilde{z}) -F_{\gamma}(\tilde{v})\text{d}y
\leq
C\int_{A_{k,t}}\left|\frac{\tilde{v}-k}{s-t}\right|^{\tilde{p}(y)}\text{d}y+C\int_{A_{k,t}}|t-s|^{\frac{\gamma
\tilde{p}(y)}{\tilde{p}(y)-\gamma}}\text{d}y,\\
&\int_{A_{k,t}}|g(\tilde{z}-\tilde{v})|\text{d}y\leq ~\int_{A_{k,t}}
\left|\frac{\tilde{v}-k}{s-t}\right|^{\tilde{p}(y)}\text{d}y+
C\int_{A_{k,t}}(|g||s-t|)^{\frac{\tilde{p}(y)}{\tilde{p}(y)-1}}\text{d}y.\notag
\end{align*}
 Therefore \eqref{5.7} becomes
\begin{align}
\int_{A_{k,s}}|\nabla \tilde{v}(y)|^{\tilde{p}(y)}\text{d}y
\leq &~C\int_{A_{k,t}\setminus A_{k,s}}|\nabla
\tilde{v}(y)|^{\tilde{p}(y)}\text{d}y+C\int_{A_{k,t}}
|\frac{\tilde{v}-k}{s-t}|^{\tilde{p}(y)}\text{d}y\notag\\
&+C\int_{A_{k,\tau}}|\tau-\sigma|^{\frac{\gamma
p_{2}}{p_{2}-\gamma}}\text{d}y+
C\int_{A_{k,\tau}}\bigg(1+(|g||\tau-\sigma|)^{\frac{p_{1}}{p_{1}-1}}\bigg)\text{d}y.\notag
\end{align}
 ``Filling the hole" and using Lemma 4.2, we obtain the desired result
\eqref{5.5}. \eqref{5.6} follows by an analogue argument.

 \emph{Second
step: Boundedness of $\tilde{v}$: estimate \eqref{5.3}}.\ \ We start by
showing that
\begin{align}
\sup\limits_{B_{\frac{1}{2}(0)}}\tilde{v}\leq
C\bigg((\int_{B_{1}(0)}(\tilde{v}^{+})^{p_{2}}\text{d}y)^{\frac{1}{p_{2}}}+1\bigg).\label{5.8}
\end{align}
 Without loss of generality we assume that $p_{1}<n$, otherwise
the assertion directly follows by the Sobolev Embedding Theorem. For
$\frac{1}{2}\leq \rho<r\leq 1$, let $\eta$ be a function of class
$C^{\infty}_{0}(B_{\frac{\rho+r}{2}})$, with $\eta\equiv 1$ on
$B_{\rho}$ and $|\nabla \eta|\leq \frac{4}{r-\rho}$. Denoting by
$p_{1}^{*}=\frac{np_{1}}{n-p_{1}}$ the Sobolev conjugate of $p_{1}$,
we introduce the quantities
\begin{align}
\varepsilon=1-\frac{p_{2}}{p_{1}^{*}}=\frac{p_{2}}{n}-\frac{p_{2}-p_{1}}{p_{1}},\
\ \beta=\varepsilon+\frac{p_{2}}{p_{1}}=1+\frac{p_{2}}{n},\ \
\theta=\varepsilon+\frac{p_{2}}{p_{1}}(1-\frac{1}{q_{-}}\frac{p_{1}}{p_{1}-1}).\notag
\end{align}
Thanks to assumption \eqref{2.7}, we have $p_{2}\leq p_{1}^{*},\
\theta>1$.

 Now we define
$
\Phi_{k,\rho}=\int_{A_{k,\rho}}(\tilde{v}-k)^{p_{2}}\text{d}y.
$
 We claim that for arbitrary $h<k$  there holds
\begin{align}
\Phi_{k,\rho}\leq &~
C\Phi_{h,r}^{\beta}\left(\frac{1}{|k-h|^{p_{2}}}\right)^{\varepsilon}
\bigg(\frac{1}{|r-\rho|^{p_{2}}}+\frac{1}{|k-h|^{p_{2}}}+\frac{1}{|k-h|^{p_{2}}}
|r-\rho|^{\frac{\gamma p_{2}}{p_{2}-\gamma
}}\bigg)^{\frac{p_{2}}{p_{1}}}\notag\\
&+C\Phi_{h,r}^{\theta}|r-\rho|^{\frac{p_{2}}{p_{2}-1}
\frac{p_{2}}{p_{1}}}\frac{1}{|k-h|^{p_{2}\theta}}.\label{5.10}
\end{align}
 Indeed, as in \cite[pp1413]{EH3}, we obtain
\begin{align}
\int_{A_{k,\rho}} (\tilde{v}-k)^{p_{2}}\text{d}y \leq
C\bigg(\int_{A_{k,\frac{\rho+r}{2}}}|\nabla
\tilde{v}|^{\tilde{p}(y)}\text{d}y+
\int_{A_{k,\frac{\rho+r}{2}}}\left|\frac{\tilde{v}-k}{r-\rho}\right|^{\tilde{p}(y)}\text{d}y
+|A_{k,r}|
\bigg)^{\frac{p_{2}}{p_{1}}}|A_{k,r}|^{\varepsilon},\label{5.11}
\end{align}
where $C=C(p_{+},p_{-})$ is a positive constant.

 Combining \eqref{5.5} and \eqref{5.11}, we derive for any $k\in\mathbb{R}$
\begin{align}
\int_{A_{k,\rho}} (\tilde{v}-k)^{p_{2}}\text{d}y \leq &~ C
|A_{k,r}|^{\varepsilon}\bigg(\int_{A_{k,r}}
\left|\frac{\tilde{v}-k}{r-\rho}\right|^{\tilde{p}(y)}\text{d}y
+\int_{A_{k,r}}|r-\rho|^{\frac{\gamma p_{2}}{p_{2}-\gamma}}\text{d}y+ \int_{A_{k,r}}(1+(|g||r-\rho|)^{\frac{p_{1}}{p_{1}-1}})
\text{d}y\bigg)^{\frac{p_{2}}{p_{1}}}
\notag\\
& +C|A_{k,r}|^{\beta}\notag\\
\leq &~C|A_{k,r}|^{\varepsilon}\bigg(\int_{A_{k,r}}
\left|\frac{\tilde{v}-k}{r-\rho}\right|^{p_{2}}\text{d}y\bigg)^{\frac{p_{2}}{p_{1}}}+C|A_{k,r}|^{\beta}+C
|A_{k,r}|^{\beta}|r-\rho|^{\frac{\gamma
p_{2}}{p_{2}-\gamma}\frac{p_{2}}{p_{1}}}\notag\\
&+C|A_{k,r}|^{\varepsilon}|r-\rho|^{\frac{
p_{2}}{p_{2}-1}\frac{p_{2}}{p_{1}}}\bigg((\int_{A_{k,r}}|g|^{q_{-}}\text{d}y)^{\frac{
p_{1}}{p_{1}-1}\frac{1}{q_{-}}}|A_{k,r}|^{1-\frac{
p_{1}}{p_{1}-1}\frac{1}{q_{-}}}\bigg)^{\frac{p_{2}}{p_{1}}}\notag\\
\leq &~C|A_{k,r}|^{\varepsilon}\bigg(\int_{A_{k,r}}
\left|\frac{\tilde{v}-k}{r-\rho}\right|^{p_{2}}\text{d}y\bigg)^{\frac{p_{2}}{p_{1}}}+C
|A_{k,r}|^{\beta}|r-\rho|^{\frac{\gamma
p_{2}}{p_{2}-\gamma}\frac{p_{2}}{p_{1}}}\notag\\
&+C|r-\rho|^{\frac{
p_{2}}{p_{2}-1}\frac{p_{2}}{p_{1}}}|A_{k,r}|^{\theta}+C|A_{k,r}|^{\beta},\label{5.12}
\end{align}
where $C=C(\lambda_{\pm},
p_{\pm},q_{-},\|g\|_{L^{q(\cdot)}(\Omega)})$ is a positive constant.

 Next,
 for $h<k$ we deduce from $u-h>k-h$ on $A_{k,r}$ that
\begin{align}
|A_{k,r}|\leq
\int_{A_{h,r}}\left|\frac{\tilde{v}-h}{k-h}\right|^{p_{2}}\text{d}y,\label{5.13}
\end{align}
and, moreover, we have
\begin{align}
\int_{A_{k,r}} (\tilde{v}-k)^{p_{2}}\text{d}y \leq \int_{A_{k,r}}
(\tilde{v}-h)^{p_{2}}\text{d}y \leq\int_{A_{h,r}}
(\tilde{v}-h)^{p_{2}}\text{d}y.\label{5.14}
\end{align}
 By \eqref{5.12}-\eqref{5.14}, we obtain
\begin{align}
\Phi_{k,\rho}\leq &~C\bigg(\int_{A_{h,r}}
\left|\frac{\tilde{v}-k}{k-h}\right|^{p_{2}}\text{d}y\bigg)^{\varepsilon}\bigg(\int_{A_{h,r}}
\left|\frac{\tilde{v}-h}{r-\rho}\right|^{p_{2}}\text{d}y\bigg)^{\frac{p_{2}}{p_{1}}}+C|r-\rho|^{\frac{\gamma
p_{2}}{p_{2}-\gamma}\frac{p_{2}}{p_{1}}}\bigg(\int_{A_{h,r}}
\left|\frac{\tilde{v}-h}{k-h}\right|^{p_{2}}\text{d}y\bigg)^{\beta}\notag\\
&+C|r-\rho|^{\frac{
p_{2}}{p_{2}-1}\frac{p_{2}}{p_{1}}}\bigg(\int_{A_{h,r}}
\left|\frac{\tilde{v}-h}{k-h}\right|^{p_{2}}\text{d}y\bigg)^{\theta}+C\bigg(\int_{A_{h,r}}
\left|\frac{\tilde{v}-h}{k-h}\right|^{p_{2}}\text{d}y\bigg)^{\beta}\notag\\
\leq &~
C\Phi_{h,r}^{\beta}\left(\frac{1}{|k-h|^{p_{2}}}\right)^{\varepsilon}
\bigg(\frac{1}{|r-\rho|^{p_{2}}}+\frac{1}{|k-h|^{p_{2}}}+\frac{1}{|k-h|^{p_{2}}}
|r-\rho|^{\frac{\gamma
p_{2}}{p_{2}-\gamma}}\bigg)^{\frac{p_{2}}{p_{1}}}\notag\\
&+C\Phi_{h,r}^{\theta}|r-\rho|^{\frac{
p_{2}}{p_{2}-1}\frac{p_{2}}{p_{1}}}\frac{1}{|k-h|^{p_{2}\theta}},\notag
\end{align}
where $C=C(\lambda_{\pm},
p_{\pm},q_{-},\|g\|_{L^{q(\cdot)}(\Omega)})$ is a positive constant.

 Our aim is now to deduce a decay estimate for the quantity
$\Phi_{k,\rho}$ to decreasing levels $k$ on balls of increasing
radii $\rho$. For this purpose we will take use of Lemma 5.1. Let us
define the sequence of levels and radii
\begin{align}
k_{i}=2d(1-2^{-i-1}),\ \ \rho_{i}=\frac{1}{2}(1+2^{-i}),\notag
\end{align}
and the quantity
\begin{align}
\vartheta_{i}=d^{-p_{2}}\Phi_{k_{i},\rho_{i}}
=d^{-p_{2}}\int_{A_{k_{i},\rho_{i}}}(\tilde{v}-k_{i})^{p_{2}}\text{d}y,\notag
\end{align}
where $d\geq 1$ is a constant that will be chosen later. First, we
note that
\begin{align}
k_{i+1}-k_{i}=\frac{d}{2}2^{-i},\ \
\rho_{i}-\rho_{i+1}=\frac{1}{4}2^{-i}.\notag
\end{align}
Exploiting \eqref{5.10} with the choice $k=k_{i+1},\ h=k_{i},\
\rho=r_{i+1},\ r=r_{i}$ and the fact that $d\geq 1$, we derive
\begin{align}
\vartheta_{i+1}=& d^{-p_{2}}\Phi_{k_{i+1},\rho_{i+1}}\notag\\
\leq &~Cd^{-p_{2}}\Phi_{k_{i},\rho_{i}}^{\beta}(d^{-1}2^{i+1})^{p_{2}\varepsilon}
\bigg((4\cdot
2^{i})^{p_{2}}+(d^{-1}2^{i+1})^{p_{2}}+(d^{-1}2^{i+1})^{p_{2}}(4^{-1}\cdot
2^{-i})^{\frac{\gamma
p_{2}}{p_{2}-\gamma}}\bigg)^{\frac{p_{2}}{p_{1}}}\notag\\
&+Cd^{-p_{2}}\Phi_{k_{i},\rho_{i}}^{\theta}(d^{-1}2^{i+1})^{p_{2}\theta}(4^{-1}\cdot
2^{-i})^{\frac{p_{2}}{p_{2}-1}\frac{p_{2}}{p_{1}}}\notag\\
=& C(d^{-p_{2}}\Phi_{k_{i},\rho_{i}})^{\beta}d^{-p_{2}(1-\beta)}
\left(\frac{d}{2}\right)^{-p_{2}\varepsilon}2^{ip_{2}\varepsilon}
\bigg(4^{p_{2}}\cdot
2^{p_{2}i}+\left(\frac{d}{2}\right)^{-p_{2}}\cdot
2^{p_{2}i}(1+(4^{-1}\cdot 2^{-i})^{\frac{\gamma
p_{2}}{p_{2}-\gamma}})\bigg)^{\frac{p_{2}}{p_{1}}}\notag\\
&+C(d^{-p_{2}}\Phi_{k_{i},\rho_{i}})^{\theta}d^{-p_{2}(1-\theta)}
\left(\frac{d}{2}\right)^{-p_{2}\theta}2^{ip_{2}\theta}(4^{-1}\cdot
2^{-i})^{\frac{ p_{2}}{p_{2}-1}\frac{p_{2}}{p_{1}}}
\notag\\
\leq &~Cd^{\frac{p_{2}}{p_{1}}(p_{2}-p_{1})}2^{ip_{2}\beta}\vartheta_{i}^{\beta}+Cd^{-p_{2}}
2^{ip_{2}\theta}\vartheta_{i}^{\theta}\notag\\
\leq &~Cd^{\frac{p_{2}}{p_{1}}
(p_{2}-p_{1})}(2^{p_{2}\beta}+2^{p_{2}\theta})^{i}
(\vartheta_{i}^{\beta}+\vartheta_{i}^{\theta}),\label{5.15}
\end{align}
where $C=C(\lambda_{\pm},
p_{\pm},q_{-},n,\|g\|_{L^{q(\cdot)}(\Omega)})$ is a positive
constant.

 Next we show that with the choice
$
d=1+\mathcal
{A}\left(\int_{B_{1}(0)}(\tilde{v}^{+})^{p_{2}}\text{d}y\right)^{\frac{1}{p_{2}}},\notag
$
where we determine the quantity $\mathcal {A}$ a bit later, the
hypotheses of Lemma 5.1 are fulfilled for the sequence
$\{\vartheta_{i}\}$.

 Due to Theorem 2.1, there exists a constant $C=C(M,p_{\pm})$ such
that
\begin{align}
d^{\frac{p_{2}}{p_{1}}(p_{2}-p_{1})}\leq C(M,p_{1},p_{2})\bigg(1+\mathcal
{A}^{\frac{p_{2}}{p_{1}}(p_{2}-p_{1})}\bigg).\notag
\end{align}
Consequently, \eqref{5.15} becomes
\begin{align}
\vartheta_{i+1} \leq c\bigg(1+\mathcal
{A}^{\frac{p_{2}}{p_{1}}(p_{2}-p_{1})}\bigg)(2^{p_{2}\beta}+2^{p_{2}\theta})^{i}
(\vartheta_{i}^{\beta}+\vartheta_{i}^{\theta}),\notag
\end{align}
where $c=c(\lambda_{\pm},
p_{\pm},q_{-},n,M,p_{1},p_{2},\Omega,\|g\|_{L^{q(\cdot)}(\Omega)})$
is a positive constant.

 On the other hand, the choice of $d$ and the
fact that $d\geq 1$ immediately yield
\begin{align}
\vartheta_{0}=d^{-p_{2}}\int_{A_{d,1}}
(\tilde{v}-d)^{p_{2}}\text{d}y\leq \mathcal
{A}^{-p_{2}}\left(\int_{B_{1}(0)}(\tilde{v}^{+})^{p_{2}}\text{d}y\right)^{-1}\int_{A_{d,1}}
(\tilde{v}-d)^{p_{2}}\text{d}y\leq\mathcal {A}^{-p_{2}}.\notag
\end{align}
 We apply Lemma 5.1 with $B=2^{p_{2}\beta}+2^{p_{2}\theta}>1,\
C=c\bigg(1+\mathcal {A}^{\frac{p_{2}}{p_{1}}(p_{2}-p_{1})}\bigg)>1,\
0<a_{1}=\theta-1<\beta-1=a_{2}$. To guarantee that the condition
$\vartheta_{0}\leq (2C)^{-\frac{1}{a_{1}}}B^{-\frac{1}{a_{1}^{2}}}$
is satisfied, we have to choose the quantity $\mathcal {A}$ in such
a way that
\begin{align}
\mathcal
{A}^{-p_{2}}=(2C)^{-\frac{1}{a_{1}}}B^{-\frac{1}{a_{1}^{2}}},\ \
\text{i.e.}\ \ \mathcal
{A}^{p_{2}(\beta-1)}=2cB^{\frac{1}{\beta-1}}\bigg(1+\mathcal
{A}^{\frac{p_{2}}{p_{1}}(p_{2}-p_{1})}\bigg).\label{5.16}
\end{align}
Note that, since
$\beta=\varepsilon+\frac{p_{2}}{p_{1}}>\frac{p_{2}}{p_{1}}$, we
always have that $p_{2}(\beta-1)>\frac{p_{2}}{p_{1}}(p_{2}-p_{1})$,
which guarantees that equation \eqref{5.16} has a unique solution
$0<\mathcal {A}\equiv \mathcal {A}(\lambda_{\pm},
p_{\pm},q_{-},n,M,p_{1},p_{2},\Omega,\|g\|_{L^{q(\cdot)}(\Omega)})<\infty$.

 In addition, we remark that global boundedness $p_{\pm}$ for $p(\cdot)$ imply
that
$p_{2}(\beta-1)=\frac{p_{2}^{2}}{n}\in[\frac{p_{-}^{2}}{n},\frac{p_{+}^{2}}{n}]$
and $\frac{p_{2}}{p_{1}}(p_{2}-p_{1})\in
[0,\frac{p_{+}}{p_{-}}(p_{+}-p_{-})]$. Furthermore, the solution
$\mathcal {A}$ of equation \eqref{5.16} depends continuously on the
parameters $p_{-}$ and $p_{+}$.

 Now Lemma 5.1 gives
$
\lim\limits_{i\rightarrow\infty} \vartheta_{i}=0,\notag
$
which, noting that $\lim\limits_{i\rightarrow\infty} \rho_{i}=\frac{1}{2}$
and $\lim\limits_{i\rightarrow\infty} k_{i}=2d$, directly translates into
$|A_{2d,\frac{1}{2}}|=0$ and therefore
$\sup\limits_{B_{\frac{1}{2}}(0)}\tilde{v}\leq 2d.$ Taking into account the
choice of $d$, we end up with
\begin{align}
\sup\limits_{B_{\frac{1}{2}}(0)}\tilde{v}\leq
C\bigg(\left(\int_{B_{1}(0)}(\tilde{v}^{+})^{p_{2}}\text{d}y\right)^{\frac{1}{p_{2}}}+1\bigg),\notag
\end{align}
where $C=C(\lambda_{\pm},
p_{\pm},q_{-},n,M,\Omega,\|g\|_{L^{q(\cdot)}(\Omega)})$.

 An argument
similar to the preceding one with the function $-\tilde{v}$, using
\eqref{5.6} instead of \eqref{5.5} yields
\begin{align}\label{5.17}
\sup\limits_{B_{\frac{1}{2}}(0)}(-\tilde{v})\leq
C\bigg(\left(\int_{B_{1}(0)}((-\tilde{v})^{+})^{p_{2}}\text{d}y\right)^{\frac{1}{p_{2}}}+1\bigg).
\end{align}
 Therefore \eqref{5.8} and \eqref{5.17} yield the desired
estimate \eqref{5.3}.

 \emph{Third step:
Boundedness of $\tilde{v}$: estimate \eqref{5.4}}.\ \ Firstly we choose
some constants we will use for our proof. By \eqref{2.7}, we know that
$\theta=\varepsilon+\frac{p_{2}}{p_{1}}\big(1-\frac{1}{q_{-}}
\frac{p_{1}}{p_{1}-1}\big)=\frac{p_{2}}{n}+1-\frac{p_{2}}{q_{-}}\frac{1}{p_{1}-1}>\frac{p_{2}}{p_{1}}$,
thus we can find a positive constant $\tilde{\alpha}$ small enough
such that
\begin{align*}
\frac{p_{2}}{p_{1}}<\frac{\theta+\tilde{\alpha}}{1+\tilde{\alpha}}\
\ \text{and}\ \
\varepsilon+\tilde{\alpha}>\tilde{\alpha}\frac{p_{2}}{p_{1}}.
\end{align*}
Then we can find positive constants
$\tilde{\beta},\tilde{\theta}$ small enough such that
\begin{align*}
\frac{p_{2}}{p_{1}}<\frac{p_{2}}{p_{1}}+\tilde{\beta}\leq
\tilde{\theta}\leq\frac{\theta+\tilde{\alpha}}{1+\tilde{\alpha}}\ \
\text{and}\ \
\beta-\tilde{\beta}-\frac{p_{2}}{p_{1}}+\tilde{\alpha}=\varepsilon-\tilde{\beta}+\tilde{\alpha}\geq\tilde{\alpha}(\frac{p_{2}}{p_{1}}+\tilde{\beta}),
\end{align*}
where the third inequality implies
\begin{align}\label{5.20}
\theta-
\tilde{\theta}+\tilde{\alpha}\geq\tilde{\alpha}\tilde{\theta}.
\end{align}
 For the above constants, it follows from \eqref{5.12}
\begin{align}
\Phi_{k,\rho}|A_{k,\rho}|^{\tilde{\alpha}} \leq &~C|A_{k,r}|^{\varepsilon}|A_{k,\rho}|^{\tilde{\alpha}}\bigg(\int_{A_{k,r}}
\left|\frac{\tilde{v}-k}{r-\rho}\right|^{p_{2}}\text{d}y\bigg)^{\frac{p_{2}}{p_{1}}}+C
|A_{k,r}|^{\beta}|A_{k,\rho}|^{\tilde{\alpha}}|r-\rho|^{\frac{\gamma
p_{2}}{p_{2}-\gamma}\frac{p_{2}}{p_{1}}}\notag\\
&+C|r-\rho|^{\frac{
p_{2}}{p_{2}-1}\frac{p_{2}}{p_{1}}}
|A_{k,r}|^{\theta}|A_{k,\rho}|^{\tilde{\alpha}}+
C|A_{k,r}|^{\beta}|A_{k,\rho}|^{\tilde{\alpha}}\notag\\
=& I_{1}+I_{2}+I_{3}+I_{4}.\label{5.21}
\end{align}
 In the following estimates we use \eqref{5.13}, \eqref{5.14}, \eqref{5.18}-\eqref{5.20}, and
the fact that $|A_{k,\rho}|\leq |A_{k,r}|\leq |A_{h,r}|$.
\begin{align}
I_{1}=& C|A_{k,r}|^{\tilde{\beta}}|A_{k,r}|^{\varepsilon-\tilde{\beta}}
|A_{k,\rho}|^{\tilde{\alpha}}\Phi_{h,r}^{\frac{p_{2}}{p_{1}}}
\bigg(\left|\frac{1}{r-\rho}\right|^{p_{2}}\bigg)^{\frac{p_{2}}{p_{1}}}
\notag\\
\leq &~C\Phi_{h,r}^{\tilde{\beta}+\frac{p_{2}}{p_{1}}}
\bigg(\left|\frac{1}{k-h}\right|^{p_{2}}\bigg)^{\tilde{\beta}}|A_{k,r}|^{\varepsilon-\tilde{\beta}}
|A_{k,r}|^{\tilde{\alpha}}\bigg(\left|\frac{1}{r-\rho}\right|^{p_{2}}\bigg)^{\frac{p_{2}}{p_{1}}}
\notag\\
\leq &~C\Phi_{h,r}^{\tilde{\beta}+\frac{p_{2}}{p_{1}}}
|A_{h,r}|^{\tilde{\alpha}(\tilde{\beta}+\frac{p_{2}}{p_{1}})}
\bigg(\left|\frac{1}{k-h}\right|^{p_{2}}\bigg)^{\tilde{\beta}}\bigg(\left|\frac{1}{r-\rho}\right|^{p_{2}}\bigg)^{\frac{p_{2}}{p_{1}}}.\label{5.22}\\
I_{2}=& C|A_{k,r}|^{\tilde{\beta}+\frac{p_{2}}{p_{1}}}|A_{k,r}|^{\beta-\tilde{\beta}-\frac{p_{2}}{p_{1}}}
|A_{k,\rho}|^{\tilde{\alpha}}|r-\rho|^{\frac{\gamma
p_{2}}{p_{2}-\gamma}\frac{p_{2}}{p_{1}}}\notag\\
\leq &~C\Phi_{h,r}^{\tilde{\beta}+\frac{p_{2}}{p_{1}}}|A_{h,r}|^{\tilde{\alpha}(\tilde{\beta}+\frac{p_{2}}{p_{1}})}
\bigg(\left|\frac{1}{k-h}\right|^{p_{2}}\bigg)^{\tilde{\beta}+\frac{p_{2}}{p_{1}}}
|r-\rho|^{\frac{\gamma p_{2}}{p_{2}-\gamma}\frac{p_{2}}{p_{1}}}.\label{5.23}\\
I_{4}\leq& C\Phi_{h,r}^{\tilde{\beta}+\frac{p_{2}}{p_{1}}}
|A_{h,r}|^{\tilde{\alpha}(\tilde{\beta}+\frac{p_{2}}{p_{1}})}
\bigg(\left|\frac{1}{k-h}\right|^{p_{2}}\bigg)^{\tilde{\beta}+\frac{p_{2}}{p_{1}}}.\label{5.24}\\
I_{3}=& C|A_{k,r}|^{\tilde{\theta}}|A_{k,r}|^{\theta-\tilde{\theta}}
|A_{k,\rho}|^{\tilde{\alpha}}|r-\rho|^{\frac{
p_{2}}{p_{2}-1}\frac{p_{2}}{p_{1}}}\notag\\
\leq &~C\Phi_{h,r}^{\tilde{\theta}}\bigg(\left|\frac{1}{k-h}\right|^{p_{2}}\bigg)^{\tilde{\theta}}
|A_{h,r}|^{\theta-\tilde{\theta}+\tilde{\alpha}}|r-\rho|^{\frac{
p_{2}}{p_{2}-1}\frac{p_{2}}{p_{1}}}\notag\\
\leq &~C\Phi_{h,r}^{\tilde{\theta}}|A_{h,r}|^{\tilde{\theta}\tilde{\alpha}}
\bigg(\left|\frac{1}{k-h}\right|^{p_{2}}\bigg)^{\tilde{\theta}}|r-\rho|^{\frac{
p_{2}}{p_{2}-1}\frac{p_{2}}{p_{1}}}.\label{5.25}
\end{align}
 Let $\tilde{\Phi}_{k,t}=\Phi_{k,t}|A_{k,t}|^{\tilde{\alpha}}$. Collecting \eqref{5.21}-\eqref{5.25}, we obtain
\begin{align}\label{5.26}
\tilde{\Phi}_{k,\rho}\leq &~C\tilde{\Phi}_{h,r}^{\tilde{\beta}+\frac{p_{2}}{p_{1}}}
\bigg(\left|\frac{1}{r-\rho}\right|^{p_{2}}+\left|\frac{1}{k-h}\right|^{p_{2}}+
\left|\frac{1}{k-h}\right|^{p_{2}}|r-\rho|^{\frac{\gamma
p_{2}}{p_{2}-\gamma}}\bigg)^{\frac{p_{2}}{p_{1}}}\bigg(\left|\frac{1}{k-h}\right|^{p_{2}}\bigg)^{\tilde{\beta}}\notag\\
&+C\tilde{\Phi}_{h,r}^{\tilde{\theta}}\bigg(\left|\frac{1}{k-h}\right|^{p_{2}}\bigg)^{\tilde{\theta}}|r-\rho|^{\frac{
p_{2}}{p_{2}-1}\frac{p_{2}}{p_{1}}},
\end{align}
where $C$ depends only on
$n,q_{-},\lambda_{\pm},p_{\pm},\|g\|_{L^{q(\cdot)}(\Omega)}$.

 To apply Lemma 5.1, taking $d\geq 1$
to be chose later and setting
$
k_{i}=d(1-2^{-i}),\ r_{i}=\frac{1}{2}(1+2^{-i}),\notag
$
we have
$
k_{i+1}-k_{i}=\frac{d}{2}2^{-i},\ r_{i+1}-r_{i}=\frac{1}{4}2^{-i}.\notag
$
 Rewriting \eqref{5.26} with $\rho=r_{i+1},r=r_{i},k=k_{i+1},h=k_{i}$ and
$\vartheta_{i}=d^{-p_{2}}\tilde{\Phi}_{k_{i},r_{i}}$ and exploiting
again the fact that $d\geq 1$, we deduce that
\begin{align}
\vartheta_{i+1}=& d^{-p_{2}}\tilde{\Phi}_{k_{i+1},\rho_{i+1}}\notag\\
\leq &~
Cd^{\frac{p_{2}}{p_{1}}(p_{2}-p_{1})}2^{ip_{2}(\tilde{\beta}+\frac{p_{2}}{p_{1}})}
\vartheta_{i}^{\tilde{\beta}+\frac{p_{2}}{p_{1}}}+Cd^{-p_{2}}
2^{ip_{2}\tilde{\theta}}\vartheta_{i}^{\tilde{\theta}}\notag\\
\leq &~Cd^{\frac{p_{2}}{p_{1}}(p_{2}-p_{1})}
2^{ip_{2}\tilde{\theta}}(\vartheta_{i}^{\tilde{\beta}+
\frac{p_{2}}{p_{1}}}+\vartheta_{i}^{\tilde{\theta}}),\notag
\end{align}
where $C$ depends only on
$n,q_{-},\lambda_{\pm},p_{\pm},\tilde{\beta},\tilde{\theta},\|g\|_{L^{q(\cdot)}(\Omega)}$.
 We now choose
$
d=1+\tilde{\mathcal
{A}}\left(\int_{A_{0,1}}\tilde{v}^{p_{2}}\text{d}y\right)^{\frac{1}{p_{2}}}|A_{0,1}|^{\frac{\tilde{\alpha}}{p_{2}}},\notag
$
where $\tilde{\mathcal {A}}$ will be fixed a bit later. Analogously
to the preceding argument we observe that
$
d^{\frac{p_{2}}{p_{1}}(p_{2}-p_{1})}\leq
c(M,p_{1},p_{2})\left(1+\tilde{\mathcal
{A}}^{\frac{p_{2}}{p_{1}}(p_{2}-p_{1})}\right).\notag
$
 The choice of $d$ gives
\begin{align}
\vartheta_{0}=d^{-p_{2}}\int_{A_{k_{0},r_{0}}}
(\tilde{v}-k_{0})^{p_{2}}\text{d}y|A_{k_{0},r_{0}}|^{\tilde{\alpha}}=d^{-p_{2}}
|A_{0,1}|^{\tilde{\alpha}}\int_{A_{0,1}}
\tilde{v}^{p_{2}}\text{d}y\leq\tilde{\mathcal {A}}^{-p_{2}}.\notag
\end{align}

 We apply Lemma 5.1 with $B=2^{p_{2}\tilde{\theta}}>1,
C=c(1+\tilde{\mathcal
{A}}^{\frac{p_{2}}{p_{1}}(p_{2}-p_{1})})>1,0<a_{1}=\tilde{\beta}+
\frac{p_{2}}{p_{1}}-1\leq\tilde{\theta}-1=a_{2}$. To guarantee that
the condition $\vartheta_{0}\leq
(2C)^{-\frac{1}{a_{1}}}B^{-\frac{1}{a_{1}^{2}}}$ is satisfied, we
have to choose the quantity $\tilde{\mathcal {A}}$ in such a way
that
\begin{align}
\tilde{\mathcal
{A}}^{-p_{2}}=(2C)^{-\frac{1}{a_{1}}}B^{-\frac{1}{a_{1}^{2}}},\ \
\text{i.e.}\ \ \tilde{\mathcal {A}}^{p_{2}(\tilde{\beta}+
\frac{p_{2}}{p_{1}}-1)}=2cB^{\frac{1}{\tilde{\beta}+
\frac{p_{2}}{p_{1}}-1}}(1+\tilde{\mathcal
{A}}^{\frac{p_{2}}{p_{1}}(p_{2}-p_{1})}).\label{5.27}
\end{align}
We note that $\tilde{\beta}>0$ which guarantees equation \eqref{5.27} has
a unique solution $0<\tilde{\mathcal {A}}<\infty$. Here
$\tilde{\mathcal {A}}\equiv\tilde{\mathcal
{A}}(n,q_{-},M,p_{\pm},\lambda_{\pm},\|g\|_{L^{q(\cdot)}(\Omega)})$.
By Lemma 5.1, we conclude that $
\lim\limits_{i\rightarrow\infty} \vartheta_{i}=0,\notag
$
which, noting that $\lim\limits_{i\rightarrow\infty} r_{i}=\frac{1}{2}$ and
$\lim\limits_{i\rightarrow\infty} k_{i}=d$, directly translates into
$|A_{d,\frac{1}{2}}|=0$ and therefore we deduce that
\begin{align}
\sup\limits_{B_{\frac{1}{2}}(0)}\tilde{v}\leq d=
C\left(\int_{A_{0,1}}(\tilde{v}^{+})^{p_{2}}\text{d}y\right)^{\frac{1}{p_{2}}}
|A_{0,1}|^{\frac{\alpha}{p_{2}}}+1,\notag
\end{align}
with $C=C(\tilde{\mathcal
{A}},n,q_{-},M,p_{\pm},\lambda_{\pm},\|g\|_{L^{q(\cdot)}(\Omega)})$.
We should note that the constant $C$ may be replaced by a constant
$C=C(n,q_{-},M,p_{\pm},\lambda_{\pm},\|g\|_{L^{q(\cdot)}(\Omega)})$.\hfill $\blacksquare$
\end{pf*}

 Now we turn to prove local boundedness for minimizers of
the functional $\mathcal {H}_{\gamma}$.
\begin{lemma}\label{Lemma 5.4}
Let $v$
 be a minimizer of $\mathcal {H}_{\gamma}$ governed by
 \eqref{5.1}. Then $v$ is locally bounded and satisfies the
 estimates
\begin{align*}
\sup\limits_{B_{\frac{R}{2}}(x_{0})} \pm v\leq
C\left(\frac{1}{|B_{R}(x_{0})|}\int_{B_{R}(x_{0})}
((\pm v)^{+})^{p_{2}}\text{d}y\right)^{\frac{1}{p_{2}}}+CR,
\end{align*}
and
\begin{align}
\sup\limits_{B_{\frac{R}{2}}(x_{0})}v\leq
C\left(\int_{B_{R}(x_{0})}((v-\kappa_{0})^{+})^{p_{2}}\text{d}y\right)^{\frac{1}{p_{2}}}
\left|\frac{A_{\kappa_{0},R}}{R^{n}}\right|^{\frac{\alpha}{p_{2}}}+R+\kappa_{0},\label{5.30}
\end{align}
for some $\alpha>0$, for all $\kappa_{0}\leq \sup\limits_{B_{R}(x_{0})}v$, where $C=C(n,L,q_{-},M,p_{\pm},\lambda_{\pm},\|g\|_{L^{q(\cdot)}(\Omega)})$.
\end{lemma}
\begin{pf*}{Proof of Lemma \ref{Lemma 5.4}}
Indeed, by the definition of $\tilde{v}$ and Lemma \ref{Lemma 5.3}, it follows
\begin{align}
\sup\limits_{x\in B_{\frac{R}{2}}(x_{0})}v(x)=& R\sup\limits_{y\in
B_{\frac{1}{2}}(0)}\tilde{v}(y)\notag\\
\leq &~
CR\bigg(\left(\int_{B_{1}(0)}(\tilde{v}^{+})^{p_{2}}\text{d}y\right)^{\frac{1}{p_{2}}}
+1\bigg)\notag\\
\leq &~
CR\bigg(\left(\int_{B_{R}(x_{0})}\left(\frac{v^{+}}{R}\right)^{p_{2}}\frac{1}{R^{n}}\text{d}x\right)^{\frac{1}{p_{2}}}
+1\bigg)\notag\\
\leq &~ C\left(\frac{1}{|B_{R}(x_{0})|}\int_{B_{R}(x_{0})}
(v^{+})^{p_{2}}\text{d}y\right)^{\frac{1}{p_{2}}}+CR.\notag
\end{align}
 Estimate \eqref{5.30} can be obtained via \eqref{5.4} by a similar argument,
taking into account that $|A_{0,R}|=R^{n}|A_{0,1}|$ and then writing
$v-\kappa_{0}$ instead of $v$.
\hfill $\blacksquare$
\end{pf*}
\begin{lemma}\label{Lemma 5.5}
Let $v$
 be a minimizer of $\mathcal {H}_{\gamma}$ governed by
 \eqref{5.1}. Then for every couple of balls $B_{\rho}\subset B_{r}\subset
 B_{R}$ having the same center $x_{0}$ and for every $k\in
 \mathbb{R}$ the following two estimates hold
\begin{align}
\int_{A_{k,\rho}}|\nabla v|^{p(x)}\text{d}x\leq
C\int_{A_{k,r}}\left|\frac{v-k}{r-\rho}\right|^{p(x)}\text{d}x
+Cr^{\lambda_{0}+n}+Cr^{n},\notag
\end{align}
and
\begin{align}
\int_{B_{k,\rho}}|\nabla v|^{p(x)}\text{d}x\leq
C\int_{B_{k,r}}\left|\frac{v-k}{r-\rho}\right|^{p(x)}\text{d}x
+Cr^{\lambda_{0}+n}+Cr^{n},\notag
\end{align}
 with $\lambda_{0}=\min\{\frac{\gamma
p_{2}}{p_{2}-\gamma},\frac{p_{1}}{p_{1}-1}(1-\frac{n}{q_{-}})\}\geq
0,\ C=C(n,L,p_{\pm},\lambda_{\pm},\|g\|_{L^{q(\cdot)}(\Omega)})$.
\end{lemma}
\begin{pf*}{Proof}
We employ an argument similar to the one used to obtain \eqref{5.5}, obtaining
\begin{align}
\int_{A_{k,\rho}}|\nabla v|^{p(x)}\text{d}x\leq &~
C\int_{A_{k,r}}\left|\frac{v-k}{r-\rho}\right|^{p(x)}\text{d}x+C\int_{A_{k,r}}|r-\rho|^{\frac{\gamma
p(x)}{p(x)-\gamma}}\text{d}x+C\int_{A_{k,r}}(|r-\rho||g|)^{\frac{
p(x)}{p(x)-1}}\text{d}x\notag\\
\leq &~C\int_{A_{k,r}}\left|\frac{v-k}{r-\rho}\right|^{p(x)}\text{d}x+C\int_{A_{k,r}}|r-\rho|^{\frac{\gamma
p_{2}}{p_{2}-\gamma}}\text{d}x\notag\\
& +C|r-\rho|^{\frac{
p_{1}}{p_{1}-1}}|A_{k,r}|^{1-\frac{1}{q_{-}}\frac{p_{1}}{p_{1}-1}}
\left(\int_{A_{k,r}}|g|^{q_{-}}\text{d}x\right)^{\frac{p_{2}}{p_{2}-1}\frac{1}{q_{-}}}+|A_{k,r}|\notag\\
\leq &~C\int_{A_{k,r}}\left|\frac{v-k}{r-\rho}\right|^{p(x)}\text{d}x+Cr^{\lambda_{0}+n}+Cr^{n}.\notag
\end{align}\hfill $\blacksquare$
\end{pf*}
 Now we shall prove H\"{o}lder regularity for the minimizers of
the functional $\mathcal {H}_{\gamma}$.
\begin{pf*}{Proof of Theorem \ref{Theorem 2.2}}
Let $v$
  be a minimizer of the functional $\mathcal {H}_{\gamma}$ governed by \eqref{5.1}. Let $ \text{osc}(v,\rho)=\sup\limits\limits_{B_\rho}v-\inf\limits_{B_\rho}v$.
 Due to Lemma \ref{Lemma 5.5}, one may proceed exactly as in \cite[Lemma 4.10]{EH3}, to see that the minimizer $v$ has also an estimate as (4.40) in \cite[Lemma 4.10]{EH3}. Again, due to Lemma \ref{Lemma 5.4} and proceeding as in \cite[Proposition 4.11]{EH3}, we have
 \begin{align}
\text{osc}(v,\rho)\leq
c\bigg(\left(\frac{\rho}{r}\right)^{\alpha_{1}}\text{osc}(v,r)+\rho^{\alpha_{1}}\bigg),\
\ \forall\  \rho<r<\frac{R}{4}.\label{5.31}
\end{align}
for some $0<\alpha_{1}<1$. By a slight modification of proof of  \cite[Proposition 4.12]{EH3},
\eqref{5.31} gives
\begin{align}
\int_{B_{\rho}}|v-(v)_{\rho}|^{p_{2}}\text{d}x\leq
C\left(\frac{\rho}{R}\right)^{n+p_{2}\alpha_{1}}\int_{B_{R}}
|v-(v)_{R}|^{p_{2}}\text{d}x+C\rho^{n+p_{2}\alpha_{1}},\notag
\end{align}
and
\begin{align}
\int_{B_{\rho}}|\nabla v|^{p(x)}\text{d}x\leq
C\left(\frac{\rho}{R}\right)^{n-p_{2}+p_{2}\alpha_{1}}\int_{B_{R}}
|\nabla v|^{p(x)}\text{d}x+C\rho^{n-p_{2}+p_{2}\alpha_{1}}.\notag
\end{align}
It follows from Lemma 5.2 that
\begin{align}
\int_{B_{\rho}}|v-(v)_{\rho}|^{p_{2}}\text{d}x\leq
C\rho^{n+p_{2}\alpha_{1}},\notag
\end{align}
and
\begin{align}
\int_{B_{\rho}}|\nabla v|^{p(x)}\text{d}x\leq
C\rho^{n-p_{2}+p_{2}\alpha_{1}}.\notag
\end{align}
Notice that each of the above inequalities combining with covering
theorem implies $v\in C^{0,\alpha_{1}}_{loc}(\Omega)$. This
concludes the proof.\hfill $\blacksquare$
\end{pf*}
\section{H\"{o}lder estimates for minimizers of functional $\mathcal {J}_{\gamma}$}
\begin{pf*}{Proof of Theorem \ref{Theorem 2.3}}
 We proceed in five
steps.

 \emph{First step: Localization}.\ \ Let
$\delta_{1}<\min\{p_{-}-1,\delta_{0}\}$ that will be chosen much
smaller a bit later. Fix a ball $B_{R_{0}}\Subset\Omega$ with
the property $\omega(8R_{0})<\frac{\delta_{1}}{4}$. Let
$B_{4R}\Subset B_{\frac{R_{0}}{4}}$. Define
$
p_{2}=\max\limits_{\overline{B}_{4R}}p(x),\
p_{1}=\min\limits_{\overline{B}_{4R}}p(x).\notag
$
We remark that by continuity of $p(x)$, there exists
$x_{0}\in\overline{B}_{4R}$, not necessarily the center, such that
$p_{2}=p(x_{0})$. Consequently we obtain
\begin{align}
&p_{2}-p_{1}\leq \omega(8R)\leq \frac{\delta_{1}}{4},\notag\\
&p_{2}(1+\frac{\delta_{1}}{4})\leq
p(x)(1+\frac{\delta_{1}}{4}+\omega(R))\leq
p(x)(1+\frac{\delta_{1}}{4}+\omega(2R))\leq p(x)(1+\delta_{1})\ \
\text{in}\ \ B_{R}(x_{0}).\notag
\end{align}
Furthermore we note the localization together with the bound \eqref{2.5}
for the modulus of continuity yields for any $8R\leq R_{0}\leq 1$:
\begin{align}
R^{-n\omega(R)}\leq \exp(nL)=c(n,L),\ \
R^{-\frac{n\omega(R)}{1+\omega(R)}}\leq c(n,L).\notag
\end{align}
In the following proofs we consider all the balls with the same
center $x_{0}$.

 \emph{Second step: Higher integrability}.\ \ By
our higher integrability result (Proposition 4.1) and localization,
it holds that
\begin{align}
\frac{1}{|B_{2R}|}\int_{B_{2R}}|\nabla
u|^{p_{2}(1+\frac{\delta_{1}}{4})}\text{d}x\leq
C_{0}\bigg(\frac{1}{|B_{2R}|}\int_{B_{2R}}|\nabla
u|^{p(x)}\text{d}x\bigg)^{1+\frac{\delta_{1}}{4}}+C_{0}
\frac{1}{|B_{2R}|}\int_{B_{2R}}\left(1+|g|^{\frac{p_{-}}{p_{-}-1}
(1+\frac{\delta_{1}}{4})}\right)\text{d}x.\notag
\end{align}
 \emph{Third step: Freezing}.\ \ Let $v\in W^{1,p_{2}}(B_{R})$ with $v-u\in W^{1,p_{2}}_{0}(B_{R})$ be
a minimizer of the functional
\begin{align}
\mathcal {G}(v)=\int_{B_{R}}f(x_{0},\nabla
v)\text{d}x=\int_{B_{R}}\tilde{h}(\nabla v)\text{d}x.\notag
\end{align}
Note that by Remark \ref{Remark 3.2} and the growth condition \eqref{2.2}, we obtain the
following estimate for the $p_{2}$ energy of $v$
\begin{align}
\int_{B_{R}}|\nabla v|^{p_{2}}\text{d}x\leq
L^{2}\int_{B_{R}}(1+|\nabla u|^{p_{2}})\text{d}x<\infty.\label{6.1}
\end{align}
Moreover, in view of \cite[Lemma 3.1]{AM}, there exist
$C=C(p_{\pm},L),\ \delta_{2}=\delta_{2}(p_{\pm},L)$ with
$0<\delta_{2}<\frac{q-p_{2}}{p_{2}}$ such that
\begin{align}
\left(\frac{1}{|B_{R}|}\int_{B_{R}}|\nabla
v|^{p_{2}(1+\delta_{2})}\text{d}x\right)^{\frac{1}{1+\delta_{2}}}\leq
C\left( \frac{1}{|B_{R}|}\int_{B_{R}}|\nabla
v|^{p_{2}}\text{d}x\right)^{\frac{1}{p_{2}}} +C
\left(\frac{1}{|B_{2R}|}\int_{B_{2R}}|\nabla
u|^{q}\text{d}x\right)^{\frac{1}{q}} ,\label{6.2}
\end{align}
for $q=p_{2}(1+\frac{\delta_{1}}{4})>p_{2}$. By the proof of Theorem \ref{Theorem 2.2}, and the boundedness of $v$, which is guaranteed by the boundedness of $u$, that there exists some $\alpha_{2}\in (0,1)$
such that
\begin{align}
\int_{B_{\rho}}|\nabla v|^{p_{2}}\text{d}x\leq
C\left(\frac{\rho}{R}\right)^{n-p_{2}+p_{2}\alpha_{2}}\int_{B_{R}}
|\nabla v|^{p_{2}}\text{d}x+C\rho^{n-p_{2}+p_{2}\alpha_{2}},\label{6.3}
\end{align}
for any $\rho$ with $2\rho<R$.

 \emph{Fourth step: Comparison
estimate}.\ \ We prove the following comparison estimate
\begin{align}
\int_{B_{R}}(\mu^{2}+|\nabla u|^{2}+|\nabla
v|^{2})^{\frac{p_{2}-2}{2}}|\nabla u-\nabla
v|^{2}\text{d}x\leq &~
C\bigg(\omega(R)\log\big(\frac{1}{R}\big)+R^{\theta_{1}}+
R^{\theta_{2}}\bigg)\int_{B_{4R}}(1+|\nabla
u|^{p_{2}})\text{d}x\notag\\
&+C\omega(R)\log\big(\frac{1}{R}\big)R^{\lambda_{1}}+CR^{\lambda_{2}}+CR^{\lambda_{3}},\label{6.4}
\end{align}
for some $0<\lambda_{1}<n,\lambda_{2}>n,\lambda_{3}>n$.

 A similar argument to the one in \cite[(4.10)]{EH2} yields
\begin{align}
\int_{B_{R}}(\tilde{h}(\nabla u)-\tilde{h}(\nabla v))\text{d}x\geq
C\int_{B_{R}}(\mu^{2}+|\nabla u|^{2}+|\nabla
v|^{2})^{\frac{p_{2}-2}{2}}|\nabla u-\nabla v|^{2}\text{d}x.\label{6.5}
\end{align}
 On the other hand, we derive
\begin{align}
\int_{B_{R}}(\tilde{h}(\nabla u)-\tilde{h}(\nabla v))\text{d}x
=& \int_{B_{R}}((f(x_{0},\nabla u)-f(x,\nabla
u))\text{d}x+\int_{B_{R}}((f(x,\nabla u)-f(x,\nabla
v))\text{d}x\notag\\
&+\int_{B_{R}}((f(x,\nabla v)-f(x_{0},\nabla
v))\text{d}x\notag\\
& =I^{(1)}+I^{(2)}+I^{(3)}.\label{6.6}
\end{align}
 We estimate $I^{(1)}$, using the continuity of the integrand with
respect to the variable $x$ (see (2.3)),
\begin{align}
I^{(1)}\leq C
\int_{B_{R}}\omega(|x-x_{0}|)\big((\mu^{2}+|\nabla
u|^{2})^{\frac{p(x)}{2}}+(\mu^{2}+|\nabla
u|^{2})^{\frac{p_{2}}{2}}\big)\big(1+|\log(\mu^{2}+|\nabla
u|^{2})|\big)\text{d}x.\notag
\end{align}
 Arguing exactly as \cite[Section 4]{EH2}, we obtain
\begin{align}
I^{(1)}\leq &~ C\omega(R) \int_{B_{R}}|\nabla
u|^{p_{2}}\log(e+\||\nabla
u|^{2}\|_{L^{1}(B_{R})})\text{d}x\notag\\
&+C\omega(R)\int_{B_{R}}|\nabla
u|^{p_{2}}\log\left(e+\frac{|\nabla u|^{p_{2}}}{\||\nabla
u|^{2}\|_{L^{1}(B_{R})}}\right)\text{d}x+ C\omega(R)R^{n}\notag\\
=& I^{(1)}_{1}+I^{(1)}_{2}+I^{(1)}_{3},\notag
\end{align}
with
\begin{align}
I^{(1)}_{1}\leq C\omega(R)\log\big(\frac{1}{R}\big)\int_{B_{R}}(1+|\nabla
u|^{p_{2}})\text{d}x.\notag
\end{align}
 Now we estimate $I^{(1)}_{2}$, using first \cite[(3.3)]{AM}, which is a
basic estimate for the $L\log L$ norm, then exploiting higher
integrability,
\begin{align}
I^{(1)}_{2}\leq &~C(p_{2},\delta)\omega(R)R^{n}\left(\frac{1}{|B_{R}|}\int_{B_{R}}|\nabla
u|^{p_{2}(1+\frac{\delta_{1}}{4})}\text{d}x\right)^{\frac{1}{1+\frac{\delta_{1}}{4}}}\notag\\
\leq &~C\omega(R)R^{n}+C\omega(R)R^{n}\left(\frac{1}{|B_{R}|}\int_{B_{R}}|\nabla
u|^{p(x)(1+\frac{\delta_{1}}{4}+\omega(R))}
\text{d}x\right)^{\frac{1}{1+\frac{\delta_{1}}{4}}}\notag\\
\leq &~C\omega(R)R^{n}+C\omega(R)R^{n}\left(\frac{1}{|B_{2R}|}\int_{B_{2R}}|\nabla
u|^{p(x)}\text{d}x\right)^{\frac{1+\frac{\delta_{1}}{4}+\omega(R)}{1+\frac{\delta_{1}}{4}}}\notag\\
&+ C\omega(R)R^{n}\left(\frac{1}{|B_{2R}|}\int_{B_{2R}}
(1+|g|^{\frac{p_{-}}{p_{-}-1}(1+\delta_{1})})
\text{d}x\right)^{\frac{1}{1+\frac{\delta_{1}}{4}}}
\notag\\
\leq &~C\omega(R)R^{n}+C\omega(R)R^{n}R^{-\frac{n\omega(R)}{1+\frac{\delta_{1}}{4}}}
\left(\frac{1}{|B_{2R}|}\int_{B_{2R}}|\nabla
u|^{p(x)}\text{d}x\right)\left(\int_{B_{2R}}|\nabla
u|^{p(x)}\text{d}x\right)^{\frac{\omega(R)}{1+\frac{\delta_{1}}{4}}}\notag\\
&+C\omega(R)R^{n}R^{-\frac{n}{1+\frac{\delta_{1}}{4}}}
\|g\|_{L^{q_{-}}(B_{2R})}^{\frac{p_{-}}{p_{-}-1}\frac{1+\delta_{1}}{1+\frac{\delta_{1}}{4}}}
R^{n[1-\frac{1}{q_{-}}\frac{p_{-}}{p_{-}-1}(1+\delta_{1})]
\frac{1}{1+\frac{\delta_{1}}{4}}}\notag\\
\leq &~C\omega(R)R^{n}+C(M)\omega(R)\int_{B_{2R}}(1+|\nabla
u|^{p_{2}})\text{d}x+C(\|g\|_{L^{q_{-}}(\Omega)})\omega(R)R^{\lambda_{1}}\notag\\
\leq &~C\omega(R)\int_{B_{2R}}(1+|\nabla
u|^{p_{2}})\text{d}x+C\omega(R)R^{\lambda_{1}},\notag
\end{align}
where
$\lambda_{1}=n-\frac{n}{1+\frac{\delta_{1}}{4}}+n[1-\frac{1}{q_{-}}
\frac{p_{1}}{p_{1}-1}(1+\delta_{1})]
\frac{1}{1+\frac{\delta_{1}}{4}}$. Notice that $\delta_{1}
<\delta_{0}<q_{-}(1-\frac{1}{p_{-}})-1$, therefore $0<\lambda_{1}<n$.

 Thus, all together we obtain
\begin{align}
I^{(1)}\leq C\omega(R)\log\big(\frac{1}{R}\big)\bigg(\int_{B_{2R}}(1+|\nabla
u|^{p_{2}})\text{d}x+R^{\lambda_{1}}\bigg).\label{6.7}
\end{align}
 We shall estimate $I^{(2)}$. By the minimizing property of $u$ and
arguing as in Section 4, we have
\begin{align}
I^{(2)}\leq &~\int_{B_{R}}(F_{\gamma}(v)-F_{\gamma}(u)+g(v-u))\text{d}x\notag\\
\leq &~
C\int_{B_{R}}|v-u|^{\gamma}\text{d}x+\int_{B_{R}}g(v-u)\text{d}x\notag\\
\leq &~ C\bigg(\int_{B_{R}}|\nabla v-\nabla
u|^{p_{2}}\text{d}x\bigg)^{\frac{\gamma}{p_{2}}}|B_{R}|^{\frac{p_{2}-\gamma}{p_{2}}}
|B_{R}|^{\frac{\gamma}{n}}+C\|g\|_{L^{\frac{p_{2}}{p_{2}-1}}(B_{R})}|B_{R}|^{\frac{1}{n}}
\bigg(\int_{B_{R}}|\nabla v-\nabla u|^{p_{2}}\text{d}x\bigg)^{\frac{1}{p_{2}}}\notag\\
\leq &~\varepsilon_{1}\int_{B_{R}}|\nabla v-\nabla
u|^{p_{2}}\text{d}x+C(\varepsilon_{1})
|B_{R}|^{(\frac{p_{2}-\gamma}{p_{2}}
+\frac{\gamma}{n})\frac{p_{2}}{p_{2}-\gamma}}+\varepsilon_{2}\int_{B_{R}}|\nabla
v-\nabla u|^{p_{2}}\text{d}x\notag\\
&+C(\varepsilon_{2})\bigg(\|g\|_{L^{\frac{p_{2}}{p_{2}-1}}(B_{R})}|B_{R}|^{\frac{1}{n}}\bigg)^{\frac{p_{2}}{p_{2}-1}}
\notag\\
\leq &~(\varepsilon_{1}+\varepsilon_{2})\int_{B_{R}}|\nabla
v-\nabla u|^{p_{2}}\text{d}x+C(\varepsilon_{1})
R^{n+\frac{p_{2}\gamma}{p_{2}-\gamma}}+C(\varepsilon_{2})
\|g\|_{L^{q_{-}}(B_{R})}^{\frac{p_{2}}{p_{2}-1}}R^{n[1+\frac{p_{2}}{p_{2}-1}(\frac{1}{n}-\frac{1}{q_{-}})]},\notag
\end{align}
where in the last but one inequality we used Young' inequality with
$
C(\varepsilon_{1})=C\left(\frac{\gamma}{\varepsilon_{1}p_{2}}\right)^{\frac{\gamma}{p_{2}-\gamma}}
\frac{p_{2}-\gamma}{p_{2}}\leq
C\left(\frac{1}{\varepsilon_{1}p_{2}}\right)^{\frac{\gamma}{p_{2}-\gamma}}
\frac{p_{2}-\gamma}{p_{2}},
$
and
$
C(\varepsilon_{2})=C\left(\frac{1}{\varepsilon_{2}p_{2}}\right)^{\frac{1}{p_{2}-1}}
\frac{p_{2}-1}{p_{2}}.
$

 Choosing $\theta_{1},\theta_{2}>0$ small enough such that
$0<\theta_{1}<p_{2}$ and
$0<\theta_{2}<np_{2}(\frac{1}{n}-\frac{1}{q_{-}})$, and setting
$\varepsilon_{i}=R^{\theta_{i}}$, we have
\begin{align}
I^{(2)} \leq &~(R^{\theta_{1}}+R^{\theta_{2}})\int_{B_{R}}|\nabla
v-\nabla u|^{p_{2}}\text{d}x+C
R^{n+\frac{p_{2}\gamma}{p_{2}-\gamma}-\frac{\gamma\theta_{1}}{p_{2}-\gamma}}
+CR^{n[1+\frac{p_{2}}{p_{2}-1}(\frac{1}{n}-\frac{1}{q_{-}})]
-\frac{\theta_{2}}{p_{2}-1}}\notag\\
\leq &~C(R^{\theta_{1}}+R^{\theta_{2}})\int_{B_{R}}(|\nabla
v|^{p_{2}}+|\nabla u|^{p_{2}})\text{d}x+CR^{\lambda_{2}}+CR^{\lambda_{3}},\notag\\
\leq &~C(R^{\theta_{1}}+R^{\theta_{2}})\int_{B_{R}}(1+|\nabla
u|^{p_{2}})\text{d}x+CR^{\lambda_{2}}+CR^{\lambda_{3}},\notag
\end{align}
where $\lambda_{2}=n+\frac{p_{2}\gamma}{p_{2}-\gamma}
-\frac{\gamma\theta_{1}}{p_{2}-\gamma}\geq
n,\lambda_{3}=n[1+\frac{p_{2}}{p_{2}-1}(\frac{1}{n}-\frac{1}{q_{-}})]
-\frac{\theta_{2}}{p_{2}-1}>n$ and in the last inequality we used
\eqref{6.1}.

 We deal with $I^{(3)}$ in a similar way to $I^{(1)}$.
Estimating in exactly the same way as in \eqref{6.7} with $v$ instead of
$u$ and doing the same splitting into $I^{(1)}$ to $I^{(3)}$, we use
higher integrability of $v$ and $u$ (\eqref{6.2} and Proposition 4.1) to
obtain
\begin{align}
I^{(3)}_{2} \leq &~C\omega(R)R^{n}\left(\frac{1}{|B_{R}|}\int_{B_{R}}|\nabla
v|^{p_{2}(1+\delta_{2})}\text{d}x\right)^{\frac{1}{1+\delta_{2}}}\notag\\
\leq &~C\omega(R)R^{n}\bigg(\left(
\frac{1}{|B_{R}|}\int_{B_{R}}|\nabla
v|^{p_{2}}\text{d}x\right)^{\frac{1}{p_{2}}} +
\left(\frac{1}{|B_{2R}|}\int_{B_{2R}}|\nabla
u|^{p_{2}(1+\frac{\delta_{1}}{4})}\text{d}x\right)^{\frac{1}{p_{2}(1+\frac{\delta_{1}}{4})}}
\bigg)\notag\\
\leq &~C\omega(R)R^{n}\bigg(\frac{1}{|B_{R}|}\int_{B_{R}}(1+|\nabla
v|^{p_{2}})\text{d}x+\left(\frac{1}{|B_{2R}|}\int_{B_{2R}}(1+|\nabla
u|^{p_{2}(1+\frac{\delta_{1}}{4})})\text{d}x\right)^{\frac{1}{1+\frac{\delta_{1}}{4}}}
\bigg)\notag\\
\leq &~C\omega(R)\int_{B_{R}}(1+|\nabla
u|^{p_{2}})\text{d}x+C\omega(R)\left(\frac{1}{|B_{2R}|}\int_{B_{2R}}|\nabla
u|^{p_{2}(1+\frac{\delta_{1}}{4})}\text{d}x\right)^{\frac{1}{1+\frac{\delta_{1}}{4}}}
\notag\\
\leq &~C\omega(R)\int_{B_{4R}}(1+|\nabla
u|^{p_{2}})\text{d}x+C\omega(R)R^{\lambda_{1}},\notag
\end{align}
where in the last inequality we used the estimate for $I_{2}^{(1)}$
to handle the second term since we assume that $B_{4R}\Subset
B_{\frac{R_{0}}{4}}$ at the beginning of the third step. All together we
end up with
\begin{align}
I^{(3)} \leq
C\omega(R)\log\big(\frac{1}{R}\big)\bigg(\int_{B_{4R}}(1+|\nabla
u|^{p_{2}})\text{d}x+R^{\lambda_{1}}\bigg).\label{6.8}
\end{align}
 From \eqref{6.5} to \eqref{6.8}, one may obtain \eqref{6.4}.

 \emph{Fifth step: Conclusion}.\ \ Now we turn to prove a decay
estimate for the $p_{2}$ energy of $u$.
We split as follows:
\begin{align}
\int_{B_{\rho}}|\nabla u|^{p_{2}}\text{d}x\leq &~
\int_{B_{\rho}}(\mu^{2}+|\nabla
u|^{2})^{\frac{p_{2}}{2}}\text{d}x\notag\\
\leq &~C\int_{B_{\rho}}(\mu^{2}+|\nabla
v|^{2})^{\frac{p_{2}}{2}}\text{d}x+C\int_{B_{\rho}}(\mu^{2}+|\nabla
u|^{2}+|\nabla v|^{2})^{\frac{p_{2}-2}{2}}|\nabla u-\nabla
v|^{2}\text{d}x\notag\\
=&  \mathcal {A}+\mathcal {B},\notag
\end{align}
where $C >0$ depends only in $p_2$.

 For $\mathcal {A}$, we deduce from \eqref{6.1} and \eqref{6.3} that
\begin{align}
\mathcal {A}\leq &~C\rho^{n}+C\int_{B_{\rho}}|\nabla
v|^{p_{2}}\text{d}x\notag\\
\leq &~CR^{n}+
C\left(\frac{\rho}{R}\right)^{n-p_{2}+p_{2}\alpha_{2}}\int_{B_{R}}
|\nabla v|^{p_{2}}\text{d}x+C\rho^{n-p_{2}+p_{2}\alpha_{2}}\notag\\
\leq &~C\left(\frac{\rho}{R}\right)^{n-p_{2}+p_{2}\alpha_{2}}\int_{B_{R}}
(1+|\nabla
u|^{p_{2}})\text{d}x+C\rho^{n-p_{2}+p_{2}\alpha_{2}}.\notag
\end{align}
 For $\mathcal {B}$, by the comparison estimate \eqref{6.4}, it follows
that
\begin{align}
\mathcal {B}\leq &~C\bigg(\omega(R)\log\big(\frac{1}{R}\big)+R^{\theta_{1}}+R^{\theta_{2}}\bigg)\int_{B_{4R}}(1+|\nabla
u|^{p_{2}})\text{d}x+C\omega(R)\log\big(\frac{1}{R}\big)R^{\lambda_{1}}+R^{\lambda_{2}}+R^{\lambda_{3}}.\notag
\end{align}
 Note that $\lambda_{1}<n<\lambda_{2},\lambda_{3}$, then we have
\begin{align}
\int_{B_{\rho}}|\nabla u|^{p_{2}}\text{d}x\leq &~
C\bigg(\left(\frac{\rho}{R}\right)^{n-p_{2}+p_{2}\alpha_{2}}+
\omega(R)\log\big(\frac{1}{R}\big)+R^{\theta_{1}}+R^{\theta_{2}}\bigg)
\int_{B_{4R}}(1+|\nabla
u|^{p_{2}})\text{d}x+CR^{\lambda_{1}}+CR^{n-p_{2}+p_{2}\alpha_{2}}.\notag
\end{align}
 On the other hand, by \eqref{2.7} we have
$n(1-\frac{1}{q_{-}}\frac{p_{1}}{p_{1}-1})>n-p_{2}$, therefore we
may choose $\delta_{1}$ and $\alpha_{2}$ small enough such that
$\lambda_{1}=n-\frac{n}{1+\frac{\delta_{1}}{4}}+n\big(1-\frac{1}{q_{-}}
\frac{p_{1}}{p_{1}-1}(1+\delta_{1})\big)
\frac{1}{1+\frac{\delta_{1}}{4}}\geq n-p_{2}+p_{2}\alpha_{2}$. Thus
\begin{align}
\int_{B_{\rho}}|\nabla u|^{p_{2}}\text{d}x\leq
C\bigg(\left(\frac{\rho}{R}\right)^{n-p_{2}+p_{2}\alpha_{2}}+
\omega(R)\log\big(\frac{1}{R}\big)+R^{\theta_{1}}+R^{\theta_{2}}\bigg)
\int_{B_{4R}}(1+|\nabla
u|^{p_{2}})\text{d}x+CR^{n-p_{2}+p_{2}\alpha_{2}}.\notag
\end{align}
 In order to apply Lemma 5.2, we may take $R_{1}>0$ small enough such that
$\omega(R)\log\big(\frac{1}{R}\big)+R^{\theta_{1}}+R^{\theta_{2}}$
smaller that $\mu$ in Lemma 5.2 for any $0<R\leq R_{1}$. Thus there
holds
\begin{align}
\int_{B_{\rho}}|\nabla u|^{p_{2}}\text{d}x\leq
C\rho^{n-p_{2}+p_{2}\alpha_{3}}\leq
C\rho^{n-p_{1}+p_{1}\alpha_{3}},\notag
\end{align}
for any $0<\alpha_{3}<\alpha_{2}$. By a standard covering argument
we deduce that $u\in \mathcal
{L}^{p_{-},n+p_{-}\alpha_{3}}_{loc}(\Omega)$, where $\mathcal
{L}^{p,\lambda}(\Omega)$ denotes Campanato's spaces, the definition
of which can be find in \cite{G}, for instance. Poincar\'{e} inequality
and a well-known property of functions in Campanato's spaces (see
\cite{G} for instance) imply that $u\in C^{0,\alpha_{3}}_{loc}(\Omega)$.\hfill $\blacksquare$
\end{pf*}
\section{$C^{1,\alpha}$ estimates for minimizers of $\mathcal {J}_{\gamma} (\gamma\in (0, 1])$}
\begin{pf*}{Proof of Theorem \ref{Theorem 2.4} $(0<\gamma\leq 1)$}
The proof
consists of three steps.

 \emph{First step: localization and freezing.}\ \ Firstly, by (2.10),
we can choose $\delta_{3}>0$ small enough such that
$\varsigma>\frac{n}{q_{-}}\frac{p_{1}}{p_{1}-1}\frac{1+\delta_{3}}
{1+\frac{\delta_{3}}{4}}$. Now let
$\delta=\min\{\delta_{0},\delta_{1},\delta_{2},\delta_{3}\}$. We
adopt the same localization argument as the proof of Theorem 2.3. In
this case all the balls $B_{CR}$ and the exponents $p_{1},p_{2}$
that we consider here are the same as in the proof of Theorem 2.3
(replace $\delta_{1}$ with $\delta$ in Section 6). Let $v\in
W^{1,p_{2}}(B_{R})$ with $v-u\in W^{1,p_{2}}_{0}(B_{R})$ be a minimizer of the functional
\begin{align}\label{7.1}
\mathcal {G}_{0}(v)=\int_{B_{R}}f(x_{0},\nabla
v)\text{d}x=\int_{B_{R}}\tilde{f}(\nabla v)\text{d}x.
\end{align}
We note that since $v$ is a minimizer of the functional $\mathcal
{G}_{0}$ with boundary data $u$ in $\partial B_{R}$, where
$u|_{\partial B_{R}}$ is the trace of a H\"{o}lder continuous
function. By Theorem 7.8 in \cite{G}, we conclude that $v\in
C^{0,\alpha_{4}}$ for some $\alpha_{4}\in(0,1)$. Therefore, for the
rest of the proof we assume that
\begin{align}
|v(x)-v(y)|\leq [v]_{\alpha_{4}}|x-y|^{\alpha_{4}}\leq
C|x-y|^{\alpha_{4}},\notag
\end{align}
holds for all $x,y\in \bar{B}_{R}$. We remark that for simplicity we
will use the same H\"{o}ler exponent for the functions $v$ and $u$,
which is not restrictive. Let us remark that, since $v$ minimizes
the functional \eqref{7.1}, by the growth condition \eqref{2.2}, higher
integrability and Remark 3.2, we obtain the following estimate for
the $p_{2}$ energy of $v$
\begin{align}\label{7.2}
\int_{B_{R}}|\nabla v|^{p_{2}}\text{d}x\leq
L^{2}\int_{B_{R}}(1+|\nabla u|^{p_{2}})\text{d}x<\infty.
\end{align}
 \emph{Second step: Comparison estimate.}\ \ We will show that
\begin{align}\label{7.3}
\int_{B_{R}}|\nabla u-\nabla v|^{p_{2}}\text{d}x\leq
CR^{\frac{\theta_{5}}{2}}\int_{B_{4R}}(1+|\nabla
u|^{p_{2}})\text{d}x,
\end{align}
for some $\theta_{5}>0$.

 Firstly we prove
\begin{align}
\mathcal {G}_{0}(u)-\mathcal {G}_{0}(v)\leq &~
C\bigg(\omega(R)\log\big(\frac{1}{R}\big)+R^{\theta_{1}}+
R^{\theta_{2}}\bigg)\int_{B_{4R}}(1+|\nabla
u|^{p_{2}})\text{d}x\notag\\
&+C\omega(R)\log\big(\frac{1}{R}\big)R^{\lambda_{1}}+CR^{\lambda_{2}}+CR^{\lambda_{3}},\label{7.4}
\end{align}
for some $0<\lambda_{1}<n,\lambda_{2}>n,\lambda_{3}>n$.

 Indeed,
since $u$ is a minimizer of the functional \eqref{1.1}, we obtain
\begin{align}
\int_{B_{R}}f(x,\nabla u)\text{d}x\leq\int_{B_{R}}f(x,\nabla
v)\text{d}x+\int_{B_{R}}(F_{\gamma}(v)-F_{\gamma}(u))\text{d}x
+\int_{B_{R}}g(v-u)\text{d}x,\notag
\end{align}
which implies
\begin{align}
\int_{B_{R}}f(x_{0},\nabla u)\text{d}x\leq &~
\int_{B_{R}}f(x_{0},\nabla v)\text{d}x+\int_{B_{R}}(f(x_{0},\nabla
u)-f(x,\nabla u))\text{d}x+\int_{B_{R}}(f(x,\nabla v)-f(x_{0},\nabla
v))\text{d}x\notag\\
&+\int_{B_{R}}(F_{\gamma}(v)-F_{\gamma}(u))\text{d}x+\int_{B_{R}}g(v-u)\text{d}x\notag\\
=& \int_{B_{R}}f(x_{0},\nabla
v)\text{d}x+I^{(4)}+I^{(5)}+I^{(6)}+I^{(7)}.\label{7.5}
\end{align}
 Arguing as $I^{(1)},I^{(3)},I^{(2)}$ in Section 6, we obtain
\begin{align}
I^{(4)}+I^{(5)}\leq
C\omega(R)\log\big(\frac{1}{R}\big)\bigg(\int_{B_{4R}}(1+|\nabla
u|^{p_{2}})\text{d}x+R^{\lambda_{1}}\bigg),\label{7.6}
\end{align}
where $0<\lambda_{1}=n-\frac{n}{q_{-}}\frac{p_{1}}{p_{1}-1}
\frac{1+\delta}{1+\frac{\delta}{4}}<n$.
\begin{align}
I^{(6)}+I^{(7)}\leq
C(R^{\theta_{1}}+R^{\theta_{2}})\int_{B_{R}}(1+|\nabla
u|^{p_{2}})\text{d}x+CR^{\lambda_{2}}+CR^{\lambda_{3}},\label{7.7}
\end{align}
where $\lambda_{2}=n+\frac{p_{2}\gamma}{p_{2}-\gamma}
-\frac{\gamma\theta_{1}}{p_{2}-\gamma}>
n,\lambda_{3}=n\big(1+\frac{p_{2}}{p_{2}-1}(\frac{1}{n}-\frac{1}{q_{-}})\big)
-\frac{\theta_{2}}{p_{2}-1}>n$.

 Therefore we may conclude \eqref{7.4} from \eqref{7.5} to \eqref{7.7}.

 Since
$\varsigma>\frac{n}{q_{-}}\frac{p_{1}}{p_{1}-1}\frac{1+\delta}{1+\frac{\delta}{4}}$,
we may choose $\theta_{3}>0$ small enough such that
$
\varsigma\geq
\frac{n}{q_{-}}\frac{p_{1}}{p_{1}-1}\frac{1+\delta}{1+\frac{\delta}{4}}+\theta_{3}.\notag
$
Again we may choose $0<\theta_{4}<\theta_{3}$ such that
\begin{align}
\varsigma+\lambda_{1}-\theta_{4}\geq n+\theta_{3}-\theta_{4}>n.\label{7.8}
\end{align}
By the assumption that $\omega(R)\leq LR^{\varsigma}$, we get
\begin{align}
\omega(R)\log\big(\frac{1}{R}\big)R^{\lambda_{1}}\leq &~
LR^{\varsigma}R^{\lambda_{1}}R^{-\theta_{4}}R^{\theta_{4}}\log\big(\frac{1}{R}\big)\notag\\
\leq &~C LR^{\varsigma+\lambda_{1}-\theta_{4}},\label{7.9}
\end{align}
for $R$ small enough.

 We deduce from \eqref{7.4}, \eqref{7.8} and \eqref{7.9} that
\begin{align}
\mathcal {G}_{0}(u)-\mathcal {G}_{0}(v)\leq
CR^{\theta_{5}}\int_{B_{4R}}(1+|\nabla u|^{p_{2}})\text{d}x,\notag
\end{align}
where $0<\theta_{5}=\theta_{5}
(\theta_{1},\theta_{2},\theta_{3},\theta_{4},
\lambda_{1},\lambda_{2},\lambda_{3},n,q_{-},p_{\pm},\varsigma,\delta),
C$ is independent of $\theta_{5}$ and $\gamma$.

 Since the grand is
of class $C^{2}$, we conclude from \cite[pp131, 137-138]{AM} that
\begin{align}
\int_{B_{R}}|\nabla u-\nabla v|^{p_{2}}\text{d}x\leq
CR^{\frac{\theta_{5}}{2}}\int_{B_{4R}}(1+|\nabla
u|^{p_{2}})\text{d}x,\label{7.10}
\end{align}
which completes the proof of \eqref{7.3}.

 \emph{Third step: Conclusion.}\ \ Firstly applying Jensen inequality we get
\begin{align}
\int_{B_{\rho}}|(\nabla u)_{\rho}-(\nabla
v)_{\rho}|^{p_{2}}\text{d}x \leq &~\int_{B_{\rho}}\left(\frac{\int_{B_{\rho}}|\nabla u-\nabla
v|\text{d}x}{|B_{\rho}|}\right)^{p_{2}}\text{d}x\notag\\
\leq &~\int_{B_{\rho}}\left(\frac{\int_{B_{\rho}}|\nabla u-\nabla
v|^{p_{2}}\text{d}x}{|B_{\rho}|}\right)\text{d}x\notag\\
\leq &~\int_{B_{\rho}}|\nabla u-\nabla v|^{p_{2}}\text{d}x.\label{7.11}
\end{align}
 Secondly, by \cite[(3.20)]{AM}, we have
\begin{align}
\frac{1}{|B_{\rho}|}\int_{B_{\rho}}|\nabla v-(\nabla
v)_{\rho}|^{p_{2}}\text{d}x\leq C\left(\frac{\rho}{R}\right)^{\beta
p_{2}} \frac{1}{|B_{R}|}\int_{B_{R}}(1+|\nabla v|^{p_{2}})\text{d}x,\label{7.12}
\end{align}
where $C>0,0<\beta<1$ and both $C$ and $\beta$ depend only on
$p_{\pm},L$.

 Now combining comparison estimate with \eqref{7.3} and
\eqref{7.11}, we deduce for any $0<\rho<\frac{R}{2}<\frac{R_{1}}{2}$
\begin{align}
\int_{B_{\rho}}|\nabla u-(\nabla
u)_{\rho}|^{p_{2}}\text{d}x\leq &~C\bigg(
\int_{B_{\rho}}|\nabla u-\nabla v|^{p_{2}}+|\nabla v-(\nabla
v)_{\rho}|^{p_{2}}+|(\nabla
v)_{\rho}-(\nabla u)_{\rho}|^{p_{2}}\text{d}x\bigg)\notag\\
\leq &~C\bigg( \int_{B_{\rho}}|\nabla u-\nabla
v|^{p_{2}}\text{d}x+\int_{B_{\rho}}|\nabla v-(\nabla
v)_{\rho}|^{p_{2}}\text{d}x\bigg)\notag\\
\leq &~C \bigg(\left(\frac{\rho}{R}\right)^{n+\beta p_{2}}
+R^{\frac{\theta_{5}}{2}}\bigg)\int_{B_{4R}}(1+|\nabla
u|^{p_{2}})\text{d}x.\label{7.13}
\end{align}
 On the other hand, we obtain (see \cite[pp133]{AM} for more details),
\begin{align}
\int_{B_{\rho}}|\nabla v|^{p_{2}}\text{d}x\leq
C\bigg(\left(\frac{\rho}{R}\right)^{n}+
\omega(R)\log\big(\frac{1}{R}\big)\bigg) \int_{B_{4R}}|\nabla
v|^{p_{2}}\text{d}x+CR^{n}.\notag
\end{align}
Therefore it follows from \eqref{7.2} that
\begin{align}
\int_{B_{\rho}}|\nabla
u|^{p_{2}}\text{d}x\leq &~C\bigg(\int_{B_{\rho}}|\nabla
v|^{p_{2}}\text{d}x+\int_{B_{R}}|\nabla v-\nabla
u|^{p_{2}}\text{d}x\bigg)\notag\\
\leq &~C\bigg(\left(\frac{\rho}{R}\right)^{n}+
\omega(R)\log\big(\frac{1}{R}\big)\bigg) \int_{B_{4R}}|\nabla
v|^{p_{2}}\text{d}x+CR^{n}+CR^{\frac{\theta_{5}}{2}}\int_{B_{4R}}(1+|\nabla
u|^{p_{2}})\text{d}x\notag\\
\leq &~C\bigg(\left(\frac{\rho}{R}\right)^{n}+
\omega(R)\log\big(\frac{1}{R}\big)+R^{\frac{\theta_{5}}{2}}\bigg)
\int_{B_{4R}}(1+|\nabla u|^{p_{2}})\text{d}x+CR^{n}.\notag
\end{align}
Thus for small $R$, applying Lemma 5.2, we obtain
\begin{align}
\int_{B_{\rho}}|\nabla u|^{p_{2}}\text{d}x\leq C\rho^{n-\tau},\notag
\end{align}
for any $\tau\in (0,1)$.

 Now let $\rho=\frac{1}{2}R^{1+\theta_{6}}$ with
$\theta_{6}=\frac{\theta_{5}}{2(n+\beta p_{2})}$, and let
$\tau=\frac{1}{4}\frac{\beta p_{2}\theta_{5}}{n+\beta p_{2}}$. Then
we deduce from \eqref{7.13} that
\begin{align*}
\int_{B_{\rho}}|\nabla u-(\nabla u)_{\rho}|^{p_{2}}\text{d}x\leq C
\rho^{\theta_{7}},
\end{align*}
with $\theta_{7}=n+\frac{\theta_{5}\beta p_{2}}{4(n+\beta
p_{2}+\frac{\theta_{5}}{2})}$. Since we can choose $\theta_{5}$
sufficient small, thus we conclude that $Du\in
C^{0,\alpha}_{loc}(\Omega)$ with
$\alpha=1-\frac{n-\theta_{7}}{p_{-}}$, which completes the proof of
Theorem 2.4 with $0<\gamma\leq 1$.\hfill $\blacksquare$
\end{pf*}
\section{Log-Lipschitz estimates for minimizers of $\mathcal
{J}_{0}$}
\begin{pf*}{Proof of Theorem \ref{Theorem 2.4} $(\gamma=0)$}
We proceed
along the lines of proof in Section 7. Notice that
$\lambda_{2}=n+\frac{p_{2}\gamma}{p_{2}-\gamma}
-\frac{\gamma\theta_{1}}{p_{2}-\gamma}= n$ with $\gamma=0$.
Therefore \eqref{7.10} becomes
\begin{align}
\int_{B_{R}}|\nabla u-\nabla v|^{p_{2}}\text{d}x\leq
CR^{\frac{\theta_{5}}{2}}\int_{B_{4R}}(1+|\nabla
u|^{p_{2}})\text{d}x+CR^{n},\notag
\end{align}
where $0<\theta_{5}'=\theta_{5}'
(\theta_{1},\theta_{2},\theta_{3},\theta_{4},
\lambda_{1},\lambda_{3},n,q_{-},p_{\pm},\varsigma,\delta), C$ is
independent of $\theta_{5}'$.

 Thus, \eqref{7.12} becomes
\begin{align}
\int_{B_{\rho}}|\nabla u-(\nabla u)_{\rho}|^{p_{2}}\text{d}x\leq & C \bigg(\left(\frac{\rho}{R}\right)^{n+\beta p_{2}}
+R^{\frac{\theta_{5}'}{2}}\bigg)\int_{B_{4R}}(1+|\nabla
u|^{p_{2}})\text{d}x+CR^{n}\notag\\
\leq & C \bigg(\left(\frac{\rho}{R}\right)^{n+\beta p_{2}}
+R^{\frac{\theta_{5}'}{2}}\bigg)\int_{B_{4R}}(1+|\nabla u-(\nabla
u)_{4R}|^{p_{2}})\text{d}x\notag\\
& \ \  +C\bigg(\left(\frac{\rho}{R}\right)^{n+\beta p_{2}}
+R^{\frac{\theta_{5}'}{2}}\bigg)\int_{B_{4R}}|(\nabla
u)_{4R}|^{p_{2}}\text{d}x+CR^{n}\notag\\
\leq & C \bigg(\left(\frac{\rho}{R}\right)^{n+\beta p_{2}}
+R^{\frac{\theta_{5}'}{2}}\bigg)\int_{B_{4R}}|\nabla u-(\nabla
u)_{4R}|^{p_{2}}\text{d}x+CR^{n},\notag
\end{align}
where $C$ depends on $M$.

 Now Lemma 5.2 implies
$
\int_{B_{\rho}}|\nabla u-(\nabla u)_{\rho}|^{p_{2}}\text{d}x\leq C
\rho^{n},\notag
$
which shows that the gradient of $u$ lies in BMO space and for any
fixed subdomain $\Omega'\Subset\Omega$, there holds
\begin{align}
\|\nabla u\|_{BMO(\Omega)}\leq
C(\Omega',n,p_{\pm},\lambda_{\pm},\|g\|_{L^{q(\cdot)}(\Omega)},M).\notag
\end{align}
 Then arguing exactly as in \cite{LdT}, one has
\begin{align}
|u(x)-u(x_{0})|\leq C|x-x_{0}|\cdot|\log|x-x_{0}||.\notag
\end{align}
The proof of Theorem 2.4 is concluded.\hfill $\blacksquare$
\end{pf*}
\begin{remark} It should be mentioned that the regularity results in \cite{EH2}, where Ekeland's variational principle was applied to the establishment of regularity in the obstacle problem associated with the functional $\int_{\Omega}f(x,u,\nabla
u)\text{d}x$, are stronger than the corresponding one in \cite{EH3}. We believe that Ekeland's variational principle can be also applied to the following heterogeneous, two-phase free boundary problem
\begin{align*}
\int_{\Omega}\big(f(x,u,\nabla
u)+F_{\gamma}(u)+gu\big)\text{d}x\rightarrow \text{min},
\end{align*}
under non-standard growth conditions, and obtain stronger regularities than the results in this paper.
\end{remark}

\end{document}